\documentclass[12pt,a4paper]{article}
\usepackage[T1]{fontenc}
\usepackage[utf8]{inputenc}
\usepackage{amsmath}
\usepackage{amsfonts}
\usepackage{amssymb}
\usepackage{setspace}
\usepackage{amsthm}
\usepackage{hyperref}
\usepackage{rotating}
\usepackage[a4paper, left=2.5cm, right=2.5cm, top=2.5cm, bottom=2.5cm]{geometry}
\usepackage{float}
\usepackage{amsmath}
\usepackage{arydshln}
\usepackage{graphicx}
\usepackage[usenames, dvipsnames]{color}
\usepackage{tabu}
\usepackage{sansmath}
\usepackage{subfigure,bm}
\usepackage[ragged]{sidecap}
\sidecaptionvpos{figure}{c}
\usepackage{natbib}
\usepackage{url}
\usepackage{currfile}
\usepackage{comment}
\usepackage{marginnote} 
                                       
\usepackage[usenames, dvipsnames]{xcolor}
\definecolor{darkblue}{HTML}{2C326E}
\usepackage{ulem} 
\newcommand\mn[1]{
}

\newcommand\labmarg[1]{
\label{#1}\mbox{} 
\mn{#1}
}

 \newcommand{\sgn}{\operatorname{sgn}} 

 \usepackage{booktabs} 

\DeclareMathOperator{\sech}{sech}
\newcommand\bone{\boldsymbol{1}}

\newcommand{\bol}[1]{\mbox{\boldmath$#1$}}

\newcommand\Ex{\text{\sf E}}
\newcommand\Prob{\text{\sf P}}

\newcommand\Cov{\text{\sf Cov}}
\newcommand\Corr{\text{\sf Corr}}
\newcommand\mL{\mathcal{L}}
\newcommand\mR{\mathbb{R}}

\newcommand{\rd}{\mathrm{d}}
\newcommand\sfleq{\sansmath{\leq}}
\newtheorem{theorem}{Theorem}
\newtheorem{proposition}{Proposition}

\theoremstyle{remark}

\newtheorem{remark}{Remark}

\usepackage{listings} 
\lstdefinestyle{mystyle}{
backgroundcolor=\color{backcolour},
commentstyle=\color{codegreen},
keywordstyle=\color{magenta},
basicstyle=\footnotesize,
breaklines=true,
keepspaces=true,
numbers=left,
numbersep=5pt,
}
\definecolor{codegreen}{rgb}{0,0.6,0}
\definecolor{backcolour}{rgb}{0.95,0.95,0.92}

\graphicspath{{figures/}} 

\begin{document}

\begin{center}
\vspace*{2cm} \noindent {\bf \large Effective computations of joint excursion times 
for stationary Gaussian processes}\\
\vspace{1cm} \noindent {\sc Georg Lindgren$^{a}$, Krzysztof Podg\'{o}rski$^{b}$,  
Igor Rychlik$^{c}$},\\
\vspace{1cm}
{\it \footnotesize  $^a$
Mathematical Statistics, Lund University}\\
{\it \footnotesize $^b$Department of Statistics, Lund
University}\\
{\it \footnotesize $^c$ Mathematical Statistics, Chalmers 
University of Technology, Göteborg} \\

\today \quad from \currfilename
\end{center}

\begin{abstract}
This work is to popularize the method of computing the distribution of the excursion times 
for a Gaussian process that involves extended and multivariate Rice's formula. 
The approach was used in numerical  implementations of the high-dimensional integration routine and in earlier work it was shown that the computations are more effective and thus more precise than those based on Rice expansions.

The joint distribution of successive excursion times is clearly related to the distribution 
of the number of level crossings, a problem that can be attacked via the Rice series expansion, 
based on the moments of the number of crossings. 
Another point of attack is the 
``Independent Interval Approximation'' (IIA) 
intensively studied for the persistency of physical systems.  
It treats  the lengths of successive crossing intervals  as 
statistically independent. Under IIA, a renewal type argument 
leads to an expression that provides the approximate interval 
distribution via its Laplace transform. 

However,  the independence is not valid in most 
typical situations. 
Even if it leads to acceptable results for the persistency exponent of the long excursion time distribution or some classes of processes, rigorous assessment of the approximation error is not readily available.
Moreover, we show that the IIA approach cannot deliver properly defined probability distributions and thus the method is limited only to persistence studies. 

The ocean science community favours a third approach, in which a class of 
parametric marginal distributions, either fitted to excursion data or 
derived from a narrow band approximation, is extended by a copula technique 
to bivariate and higher order distributions.   

This paper presents an alternative approach that is both more general, more accurate 
and relatively unknown. It is based on exact expressions for the probability density for 
one and for two successive excursion lengths. The numerical routine {\sf RIND} computes 
the densities using recent advances in scientific computing and is easily accessible for 
a general covariance function, via simple {\sc Matlab} interface. 

The result solves the problem of two step excursion dependence for a general stationary 
differentiable Gaussian process, both in theoretical sense and in practical numerical sense.    
The work offers also some analytical results that explain the effectiveness 
of the implemented method.
\end{abstract}

\noindent ASM Classification:  60G15, 60G10, 58J65, 60G55, 62P35, 65C50, 65D30

\section{Introduction}
\subsection{The problem and some of its early history}
One of the central problems in Rice's second article on random noise, \citep{Rice_b}, 
is the statistical characterization of the zeros of a stationary Gaussian process.  
Rice's formula for the expected number of zeros, and more generally, 
of non-zero level crossings, is a first step, but Rice also presents the in- and exclusion series, 
the ``Rice series'', for the distribution of the time between two successive 
mean level crossings. 

The distribution of the number of zero crossings is naturally connected to the distribution 
of the time series of successive lengths of excursions above and below zero.  
Both lines of approach were followed during the decades following Rice's article. 

\citet{LonguetHiggins1962,LonguetHiggins1963} improved considerably on the original 
Rice series for the number of crossings and derived a rapidly converging moment series 
for the probability density of zero crossing intervals. He also compared approximations 
based on the initial terms in the series, with experimental results, 
\citep{FavreauEtAl1956}, and with earlier alternative series, suggested by 
\citet{McFadden1956,McFadden1958}. Early experiments with the series of 
zero crossing intervals were also performed by \citet{Blotekjaer1958}. 

The studies by \citet{McFadden1958} and \citet{Rainal1962}, with more details 
in \citep{Rainal1963}, are of particular interest for the present article, since they contain 
systematic theoretical as well as experimental studies of the dependence between 
successive crossing intervals. Three approximation candidates were studied, 
independence, ``quasi''-independence, which i.a.\ assumes that the sum of two 
successive intervals is independent of the next one, and Markov dependence. 
The first two cases were analysed by renewal type arguments and Laplace transforms, 
\citep{Cox1962,Sire2008},
and numerical solutions were compared to experiments. 
The Markov assumption, first suggested by \citet{McFadden1957AMS}, was tested by 
variance and correlation parameters against experiment. All three assumptions were 
rejected for general Gaussian processes. 

The tradition with experimental testing of the dependence assumptions, including the Markov 
assumption, was continued by \citet{Mimaki1973} and co-workers, \citep{MimakiTanabeWolf1981,MimakiSatoTanabe1984,MimakiMyokenKawabata1985}. 
\citet{Munakata1996} listed solved and unsolved crossing problems, focusing on 
experimental evidence and practical application of the available traditional methods 
to noise in signals. 

On the theoretical side, \citet{CramerLeadbetter1967} derived formulas for crossing 
moments of arbitrary order under minimal assumptions for Gaussian processes, and \citet{Zahle} gave 
a generalized Rice's formula for non-Gaussian processes. 
\citet{Lindgren1972AAP} introduced a regression technique, 
\citep{Slepian1963}, for the excursion length in a differentiable Gaussian process, a 
technique that formed a first step towards the numerical algorithms that will be used 
in this paper.     

\subsection{Renewed interest and new exact tools}
During the years around 1990 the interest in crossing interval distributions and their tail behaviour 
increased in material science, optics, statistical physics, and other areas, \citep{Brainina2013}. 
In this work, the emphasis was on the tail distribution of the crossing interval often 
referred to as {\sl persistency} and, in particular, on its the rate of the convergence to zero, 
as expressed by the {\sl persistency exponent}.
The ``independent interval assumption'' (IIA) was applied both to Gaussian processes and to 
other process models, and compared to experiments. \citet{Sire2008} describes the renewal 
and Laplace transform arguments, and gives many references from the physics literature. 
The diffusion processes and their persistency exponent was analyzed in \citep{MajumdarSBC}.
In the following development, involving experiments, 
simulations, and theoretical analysis of diffusion phenomena has led to deepened insight into  
the asymptotic properties of crossing distribution for a range of stochastic processes, 
as conveniently surveyed in
\citep{BrayMajumdarSchehrAiP2013}, and with recent advances given in 
\citep{PoplavskyiSchehr2018}, where an important explicit form of the persistency coefficient 
for the diffusion of order 2 has been obtained by rather deep combination across different 
developments in theoretical physics.

In other than physics areas of research, one should notably point to the ocean 
science and engineering literature. There in the studies of metaocean, the time on successive crossing periods has been studied for the particular spectra occurring in different sea states
\citep{WistMyrhaugRue2004}.
Parametric spectra has been proposed tying the natural condition at geographical location and at the sea state at given time and time crossing distribution has been elaborated in many examples \citep{Ochi1998}. 

During the same years, new tools were developed in applied probability and in scientific computing. 
\citet{Durbin1985JAP} gave the exact expression for the first passage density of a 
non-differentiable Gaussian process to a general level. The result was generalized to smooth 
processes by \citet{Rychlik},  who also expanded the formula to give the exact probability 
density of  excursion intervals for a differentiable Gaussian process, 
\citep{Rychlik1990SPA}. The exact formula gave the marginal distribution only, but 
\citet{Rychlik1987JAP} also used the regression technique, described in \citep{LindgrenRychlik1991ISR}, 
to derive an almost correct density for the joint distribution of two successive intervals. 

\citet{PodgorskiRM} presented the exact formula for the joint 
density of two or more successive crossing intervals. 
The formulas, to be described in Section~\ref{s:exact}, involve the conditional 
expectations of the derivatives at crossings and the indicator functions that the process 
stays above or below the level in the intervals between crossings. No analytic expressions for these 
expectations are known, but they are readily computable by the high-dimensional integration methods 
that have been made possible by advances in scientific computing. 

By proper use of numerical linear algebra and numerical integration techniques \citet{Genz1992} 
and \citet{Rychlik1992CSS} almost simultaneously developed practically useful routines for 
computation of high-dimensional normal integrals. Genz's routines were expanded to 
very high dimensions, \citep{GenzKwong2000}, while \citet{PodgorskiRM}, amended 
Rychlik's routine to include conditioning on level crossings and derivatives. 
The routine, called {\sf RIND}, was included in the {\sc Matlab} toolbox {\sf WAFO}; see 
\citep{BrodtkorbJLRRS} and \citep{gitwafo}.
\citet{Brodtkorb2006} combined all the described ideas into a powerful  
computational tool, adding new tests to control accuracy, and embedding it in user-friendly code 
for use on Gaussian process crossing problems.  
\section{A background on interval dependence and the IIA}\label{sec:IIA}
\subsection{A smooth process and its clipped version}\label{subsec:clipped}
\subsubsection{Distribution relations}\label{sss:distrel}
\begin{figure}[tbh]
\centerline{
\includegraphics[width=0.95\textwidth]{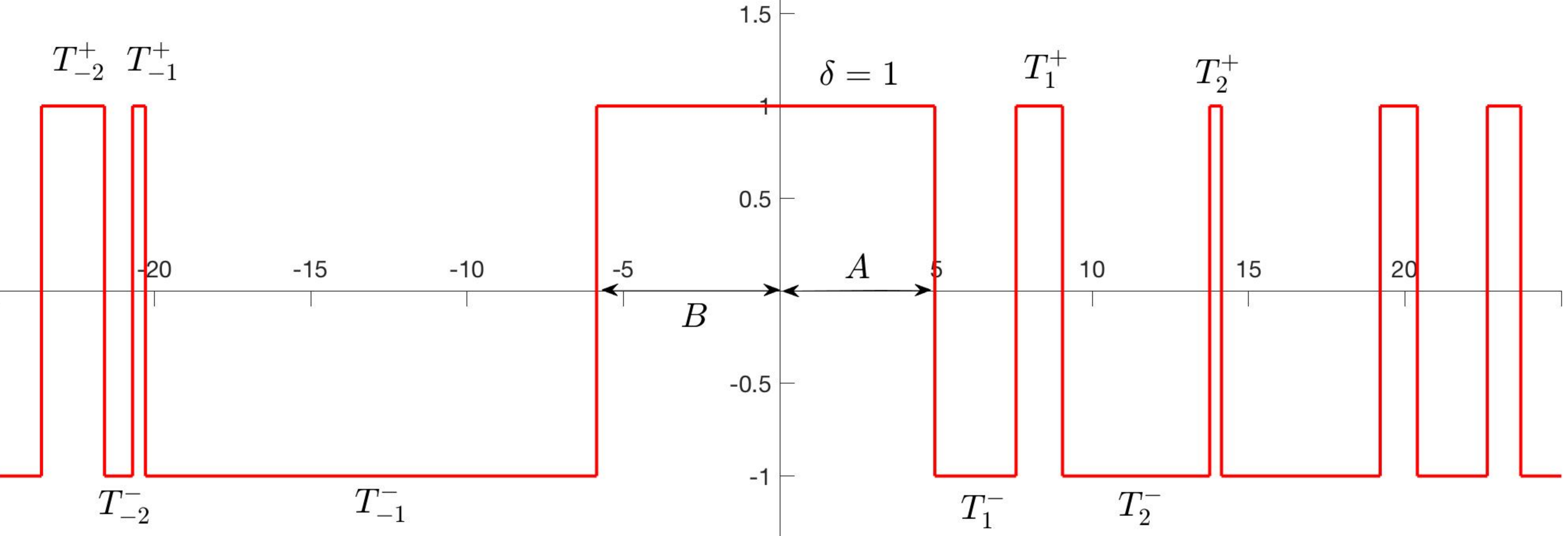}
}
\caption{ A clipped process $D_c$ with the interval $T_0^{\pm}$ split in backward and forward delays 
$B$ and $A$, respectively. The case of the origin state $\delta=1$.}
\label{DelSw}
\end{figure}

We consider a smooth process $X(t), -\infty < t < \infty$, and the time instants of $u$-level crossings, 
$S_i$, leading to two sequences of interlaced intervals of lengths  
$T_i^+$, $i=\pm 1, \pm 2, \dots$,  for the excursions above $u$-level and $T_i^-$, 
$i=\pm 1, \pm 2, \dots$, for the analogous excursions below $u$. 
We label the interval that contains the origin $T_0^{\pm} = A+B$; it may be an excursion above or below $u$. 
Figure~\ref{DelSw} explains the principle for indexing. The variable $\delta$ is introduced to keep track 
of excursions above, $\delta = 1$, or below $u$, $\delta = -1$. The process $D_c(t) = +1(-1)$ when 
$X(t)>u(<u)$ is called the {\it clipped} version of the $X$-process at level $u$. The clipped process is a 
special case of a {\it switch} process that switches between states $+1$ and $-1$ at random times points.  

If the smooth process $X(t)$ is stationary, i.e.\ its distribution is unchanged after a shift of time, 
its clipped version $D_c(t)$ is also stationary. The {\it point process} of $u$-level crossings $\{S_i\}$,  is a 
stationary point process, the joint distribution of the number of points in disjoint time intervals 
only depends on the length and relative locations of the time intervals, 
not on their absolute locations. 

We now turn to the distributions of the interval lengths $T_i^+$ and $T_i^-$, 
which are uniquely determined by the distribution of the $X$-process. This needs some care. 
In Figure~\ref{DelSw} we split the interval that contains the origin in two parts, with a forward delay time 
$A$ to the first crossing on the positive side, and a backward delay time $B$ since the last crossing 
on the negative side. 

We seek the relation between the distribution of $A, B$ and the interval distributions 
observed in an infinitely long realization of an ergodic process. 
Denote, for a fixed level $u$, by $f_T^+(t)$ and $f_T^-(t)$ the probability densities 
of the excursions above and below the level, respectively.  For a Gaussian process, 
when $u$ is equal to the mean level, the two densities are equal, 
$f_T(t) = f_T^+(t) = f_T^-(t)$ and the clipped process is symmetric with respect to the abscissa.  
Let $\mu^+, \mu^-$, and $\mu$ denote the mean interval lengths in the asymmetric and symmetric cases. 


For a symmetrically clipped process we introduce the following densities: 
 $f_{A,B}(t)$, the joint density of the forward delay $A$ and the backward delay $B$; 
$f_{A+B}(t)$, the density of the interval that contains the origin;  
$f_A(t)=f_B(t)$, the marginal densities of the forward and backward delay times.
The simple expressions are
\begin{subequations}\label{A}
\begin{align}
f_{A,B}(a,b) &= \frac{f_T(a+b)}{\mu}, \label{Aa} \\
f_{A+B}(t) &= \frac{t f_T(t)}{\mu},  \label{Ab} \\
f_A(t) = f_B(t) &= \frac{\int_t^\infty f_T(s) \, \rd s}{\mu}, \label{Ac}
\end{align}
\end{subequations}
with the intuitive interpretation of \eqref{Aa} that the location of the origin relative to the endpoints of the 
``center'' interval is uniform; see \citet[Exercise~13.3.2]{DaleyVereJones2008}. and Appendix~\ref{app:a3}.
Obviously we can conclude from \eqref{Ac} that if inter-crossing time 
and first crossing time 
have the same distribution, then this distribution is necessarily exponential. The converse is evident. 

For the asymmetric case, with $\delta = \pm 1$ indicating the status of the interval that contains 
the origin, the distribution of $(A,B,\delta)$ is given through
\begin{equation}
\begin{split}
\Prob{(\delta = 1)}=\frac{\mu^+}{\mu^-+\mu^+}, &~
\Prob{(\delta = -1)}=\frac{\mu^-}{\mu^-+\mu^+}, \\[0.3em]
f_{A,B|\delta}(a,b \mid 1)=\frac{f^+(a+b) }{\mu^+}, &~ 
f_{A,B|\delta}(a,b \mid -1)=\frac{f^-(a+b) }{\mu^-}, 
\end{split}\label{B}
\end{equation}
where $f_{A,B|\delta}$ stands for the conditional density.

The discrepancy between the long run interval distributions in a stationary point process and 
distributions taken from a frozen starting point has been first discussed for telephone calls and 
solutions have been worked out in the context of Palm measures, renewal processes, 
horizontal window conditioning, and the Rice formula. 
The mathematical foundations have been long resolved, see \citep{Palm}, \citep{Khinchin},  
\citep{RyllNardzewski}, and \citep{DaleyVereJones2008}, as well as in the key renewal theorem.

\begin{figure}
\begin{center}
\includegraphics[height=80mm]{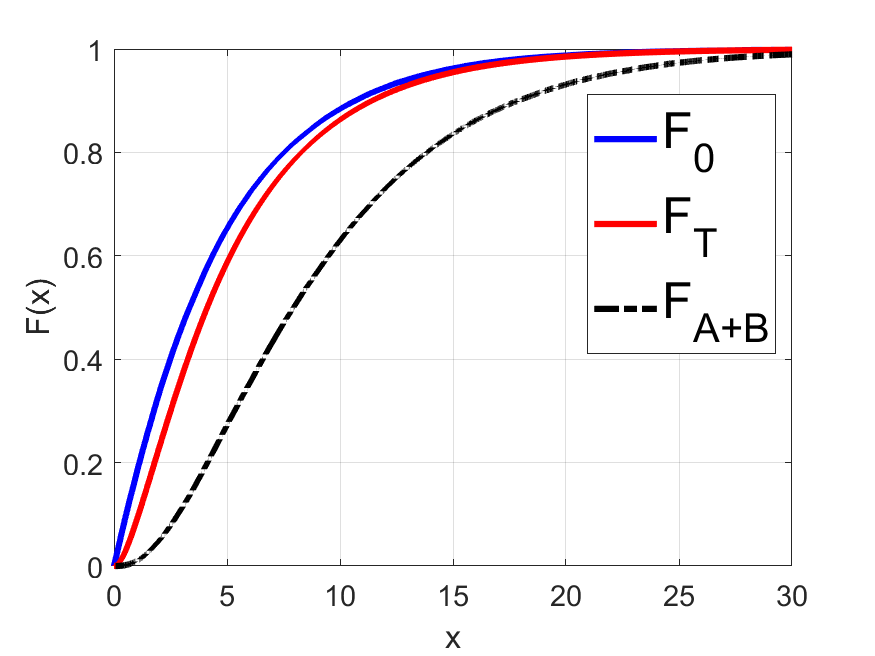}
\end{center} 
\caption{
Cumulative distribution functions (cdf) for the distance from the origin to the 
first crossing (blue), in-between switch times including the origin (black,  dash-dotted), and   
all in-between switch times (red) for a Gaussian switch process {\sf LH1}, in Table~\ref{Tab:1}. }
\label{fig:inspection}
\end{figure}

The clipped process is stationary by its design. However, we observe that for the intervals 
of its constant values (plus or minus one)  that include the origin of the horizontal axis are 
not distributed the same as the intervals of constant values (excursion intervals) as observed 
over the entire real line. 
This can be observed in Figure~\ref{fig:inspection}, where the distribution of the excursion 
intervals including the origin are contrasted with the excursion intervals as collected over the whole line. 
In fact, as observed in the graphs, the distribution of the in-between intervals 
over the whole line is closer to the distribution of the distance from the origin to the first 
crossing rather than the distribution of the entire in-between interval containing the origin. 
We conclude that statistically speaking the origin of the horizontal line hits larger intervals than those 
following from the distribution of the excursion times in agreement with the well known inspection 
paradox.
\subsubsection{Covariance function and its Laplace transform}

For a smooth stationary process  $X(t)$, we define a clipped process at the level $u$ as a process 
$D_c(t)$ that takes  value one when $X(t)>u$ and value minus one when $X(t)<u$. 
It is obvious that $D_c$ is also a stationary process and its intervals of constant values  
contains information about the length of excursions above the level $u$. 
Therefore properties of the clipped process have been used for analysis of  persistence of 
the underlying process $X(t)$.

It is rather obvious that the covariance of the clipped process can be written as follows 
\begin{align}
R_u^c(t)&=\Cov{(D_c(t+s), D_c(s))}
=\Prob{(X(t)>u,X(0)>u)} +\Prob{(X(t)<u,X(0)<u)} \nonumber \\
&~\hspace{5mm}-\Prob{(X(t)<u,X(0)>u)}-\Prob{(X(t)>u,X(0)<u)}- \left(1-2F_X(u)\right)^2 \nonumber \\
&=4F_X(u)\left(\Prob{(X(t)<u \mid X(0)<u)}- F_X(u)\right), \label{ClippedCovariance}
\end{align}
where $F_X$ is the cdf of $X(0)$. 

Let us additionally assume that the process $X(t)$ is Gaussian with covariance $R_X(t)$ 
and zero mean. In this case, if we denote $\widetilde u=u/\sqrt{R_X(0)}$ and 
$\rho_t=R_X(t)/R_X(0)$, then the respective probabilities can be written as 
\begin{align*}
\Prob{(X(t)<u \mid X(0)<u)}&=\Prob(\rho_tZ+\sqrt{1-\rho_t^2}\,Y<\widetilde u \mid Z<\widetilde u)\\
&=\Ex{\left({\Phi}\left(\frac{\widetilde u- \rho_tZ}{\sqrt{1-\rho_t^2}}\right)\, \Big | \: Z<\widetilde u\right)},
\end{align*}
where $Z$ and $Y$ are independent standard normal variables. 
In the Gaussian case 
\begin{align*}
R_u^c(t)&=4\Phi(\widetilde u)
\left(\Ex{\left({\Phi}\left(\frac{\widetilde u- \rho_tZ}{\sqrt{1-\rho_t^2}}\right)
\, \Big | \, Z<\widetilde u\right)}- \Phi(\widetilde u)\right)\\
&=4\Phi(\widetilde u)\left(\frac{1}{\sqrt{2\pi} 
\Phi(\widetilde u)}\int_{-\infty}^{\widetilde u}e^{-\frac{z^2}{2}} {\Phi}\left(\frac{\widetilde u- \rho_tz}{\sqrt{1-\rho_t^2}}\right) \, \rd z - \Phi(\widetilde u)\right),
\end{align*}
where $Z$ is a standard normal random variable. 

The formula for the auto-covariance takes a particularly simple explicit form for the symmetric case of level $u=0$, $R_0^c(t)=\frac 2 \pi \arcsin \rho_t$. 
Since the covariance $\rho_t$ defines a Gaussian process up to a scaling constant, we see, somewhat surprisingly, that clipping a Gaussian process does not lose any structural information about the original process.  

To apply IIA we will link the structure of the clipped process to that of 
a switch process with {\sl independent intervals}. The link will be the covariance functions, or rather their Laplace transforms;  for the clipped process, 
$\mL R_u^c (s) = \int_0^\infty e^{-st} R_u^c(t) \, \rd t$, and for the Gaussian symmetric case,
\begin{align}
\mL \left( \frac{2}{\pi} \arcsin \rho_t \right), \label{LRG}
\end{align}
that will be matched with the Laplace transform $\mL R $ of the stationary covariance $R(t)$ of a switch process with independent intervals.
\subsection{A renewal point process and its switch process}\label{subsec:switch}
\subsubsection{A switch process}
The clipped process $D_c(t)$ is generated from a smooth stationary process $X(t)$ via the 
point process of $u$-level crossings. That mechanism imposes certain restrictions on its 
distribution. For example, it is well-known that no non-trivial Gaussian process can have exactly independent 
mean level crossing intervals, see \cite{palmer_1956,McFadden1958,LonguetHiggins1962}.  

A start from a general simple stationary point process gives more flexibility for a switch $\pm$~process. 
To distinguish the construction from the clipping procedure we consider a stationary 
{\it marked} point process $\{S_i\}$ on the real line where a sequence of alternating marks 
$\varepsilon_i$ attached to the points indicate if the switch is $-/+$ or $+/-$. 
The distances between a $-/+$ switch at $S_i$ and the following $+/-$ switch is 
labelled $T_i^+ = S_{i+1} - S_i$, while the next switch interval, from $+/-$ to $-/+$, 
is denoted $T_i^-$. We denote the switching process by $D_s(t)$, and introduce the 
conditional probabilities 
\begin{align*}
P_\delta(t) &= \Prob (D_s(t) = 1 \mid D_s(0) = \delta), \; \delta = \pm 1.
\end{align*}
Note that relations \eqref{A} and \eqref{B} hold for $D_s(t)$
\subsubsection{A renewal switch process and its covariance function}
If all $T_i^+$ and $T_i^-$ are independent we have an alternating (delayed) renewal process, and if  
the distribution of the centre interval is given by \eqref{B} then we have a stationary switching process, 
and it has the ``Independent Interval Property'', IIP, (as different from Approximation). 
Such a process has covariance function $R^s(t)=\Cov{(D_s(t_0), D_s(t_0+t))}$ of the form
\begin{align*}
R(t) &= \frac{2}{\mu^++\mu^-}
 \left(
 P_{1}(t)\mu^+-P_{-1}(t)\mu^-+
\mu^+ 
 \frac{\mu^--\mu^+}{\mu^++\mu^-}
 \right),
 \end{align*}
The Laplace transform of the covariance is given by the Laplace transforms of 
the interval distributions, $\Psi_{\pm} (s) = \int_0^\infty e^{-st} f_{\pm}(t) \, \rd t $ and the mean interval 
lengths, $\mu_+$ and $\mu_-$,
\begin{align}
\mL R(s) &= \frac{4}{s\left(\mu^++\mu^-\right)}\left(\frac{\mu_+\mu_-}{\mu_++\mu_-} - 
\frac1s \frac{(1-\Psi_+(s))(1-\Psi_-(s))}{1-\Psi_-(s)\Psi_+(s)} \right). \label{LS}
\end{align}
In the case when the distributions of $T_i^+$ and $T_i^-$ are the same, we obtain the two relations
\begin{align}
\mL R(s) &= \frac{1}{s}\left (1 - \frac{2}{s\mu}\: \frac{1 - \Psi (s)}{1 + \Psi(s)}\right),\label{LSsym} \\[0.3em]
\Psi(s) &=\frac{2-s\mu(1-s \, \mL R(s))}{2+s\mu(1-s \, \mL R(s))}. \label{215}
\end{align}
The above result is in agreement with formula (215) in \citet{BrayMajumdarSchehrAiP2013}.

\subsection{The persistence exponent and the IIA}

\subsubsection{The IIA principle}
We can now formally state the  fundamentals of the IIA approach. \vspace{2mm}

\noindent {\sc The IIA principle for symmetric crossing distance: }{\sl 
Find the covariance function $R_c(t)$ \eqref{ClippedCovariance} of 
the clipped process from its distribution 
and match its Laplace transform to that of a symmetric renewal process \eqref{LSsym},
\begin{align}
\label{eq:cov_match}
\mL R_c (s) &= 
\mL R(s) = \frac{1}{s}\left (1 - \frac{2}{s\mu}\: \frac{1 - \Psi (s)}{1 + \Psi(s)}\right ),
\end{align}
and solve for $\Psi(s)$, according to  \eqref{215}. 
If the clipped process is Gaussian, set $\mL R_c(s) = \mL \left( \frac{2}{\pi} \arcsin \rho_t \right)(s)$. 
 \vspace{2mm}
}

We note that in the symmetric case the inter-switch distribution is a function only of the Laplace 
transform of the covariance function. 
For the asymmetric case, \eqref{LS}, there are two distributions, 
$f_+$ and $f_-$, to be matched to the covariance function.  
Thus if we do have the covariances of the switch process, we need one more relation 
to solve for these distributions. 
For that different strategies could be taken. For example, one can match two Slepian models for upcrossing and downcrossing 
with the non-stationary mean of the non-delayed switch process. 
In the result one could determine both $f_+$ and $f_-$, 
see also \citep{Sire2007} and \citep{Sire2008} for the related approach.

When closely examined, the approach is somewhat mathematically 
inconsistent as it uses the poles of a Laplace transform of a function that is not a probability 
distribution  to approximate the distributional tails  of the inter-crossing intervals. 
This is because the solution $\Phi$ of \eqref{eq:cov_match} does not, in general, correspond to a valid probability distribution. 
In Appendix~\ref{Val-IIA}, we further elaborate on these issues.

Nevertheless, the IIA principle, while not delivering a valid approximation of the distribution, works fairly well in approximating its tail behavior. 
The latter is best described in the terms of the {\sl persistency exponent}. 

\subsubsection{The persistence exponent}
The most common meaning of persistence is as the tail of the first-crossing distribution
\begin{align}
Q_T &= \Prob (X(t) \text{ does not change sign between $t=0$ and $t=T$}),
\label{QT}
\end{align} 
in particular for large $T$. Alternatively, the probability that the process stays above (or below) the level $u$
in the entire interval $[0, T]$.  

For Gaussian processes the asymptotic behaviour of the persistence depends only on the covariance 
function. General results about the asymptotic persistence are scattered and not very precise. 
The most precise statements about the decay have been formulated for processes with 
non-negative correlation function, while for oscillating correlation only upper and lower bounds have been 
obtained.

\paragraph{Non-negative correlation:} 
\cite[Thm.~1.6]{DemboMukherjee2015} give a precise meaning to the ``exponential tail'' property: 
If the correlation function $r(t)$ of a stationary centered Gaussian process $X(t)$ is everywhere non-negative, 
then there exists a non-negative limit 
\begin{align}
b_r &= - \lim_{T\to \infty} \frac1T \log \Prob\left(\inf_{t\in [0,T]} X(t) >0\right).
\end{align}
\cite{FeldheimFeldheim2015} comments that $b_r$ is necessarily finite, and thus 
it is meaningful to formulate the persistence tail as 
\begin{align}
Q_T &= \Prob\left(\inf_{t\in [0,T]} X(t) >0\right)+\Prob\left(\sup_{t\in [0,T]} X(t) <0\right)= 2 e^{-(\theta + o(T))T} = e^{-(\theta + o(T))T}, 
\label{exponentialindex}
\end{align}
where $o(T)\to 0$ as $T$ increases without bound (Landau's little o symbol), with $0 \leq \theta = b_r$ as the {\sl persistence exponent}. 
One should note that since nothing is said about 
the asymptotics of $o(T)\cdot T$, relation \eqref{exponentialindex} only gives the main order of decay 
in the sense that for every $\varepsilon > 0$, $Q_T / \exp(-(\theta \pm \varepsilon)T) \to (0, \infty)$.  

\paragraph{Oscillating correlation:}  Very little is known about the persistence for oscillating correlation.  
\cite{AntezanaEtAl2012} studied the low-frequency white noise process with correlation function 
$r(t)=\sin(t)/t$ and proved the existence of exponential upper and lower bounds,
\begin{align}
0 < e^{-cT} &\leq Q_T \leq e^{-CT}, \text{ with $c, C > 0$}. \label{exponentialbounds}
\end{align}
\cite{FeldheimFeldheim2015} generalized \eqref{exponentialbounds} to processes whose spectral measure 
is bounded away from zero and infinity near the origin. The condition is automatically fulfilled if the process 
has a spectral density  $S(\omega)$ with $m < S(\omega) < M$ for all $\omega \in [-a, a]$ for some finite 
$a, m, M >0$. 
 
\paragraph{Persistence approximation via IIA:} 
The IIA attempts to approximate the inter-crossing distance distribution through its Laplace transform 
and the inverse transform. 
While it fails to deliver the proper distribution, it still yields decent approximations  of the tail behavior.  
In fact, the inter-crossing distance is often close to an exponential
density $f_T(t) = \theta e^{-\theta t}, t>0$, with Laplace transform 
$\Psi (s) = \frac{\theta}{\theta+s}$. For the exponential distribution $Q_T = e^{-\theta T}$ 
an approximate value 
$\theta_{IIA}$ for $\theta$ is minus the largest pole 
of $\Psi(s)$ in \eqref{215}, i.e.
\begin{align}
\theta_{IIA} = - \max \{s; 2+s\mu(1-s \, \mL R(s)) = 0\},\label{ThetaIIA}
\end{align}
\cite[Eqn.~(217)]{BrayMajumdarSchehrAiP2013} for the Gaussian case:

As noted in Section~\ref{sss:distrel}, exponential inter-crossing distance implies exponential first crossing 
distance with the same parameter. Thus, the persistence exponent $\theta_{per}$, defined by 
$\Prob ( X(t) \text{ does not change sign between $t=0$ and $t=T$}) \propto e^{-\theta_{per} T}$ 
for large $T$, can be approximated as $\theta_{per} = \theta_{IIA}$. 
In Appendix~\ref{Val-IIA}, we present the IIA 
approximation of the persistency exponent for  the diffusion in the dimension two.
\begin{figure}[t]
\centerline{
\includegraphics[width=0.9\textwidth]{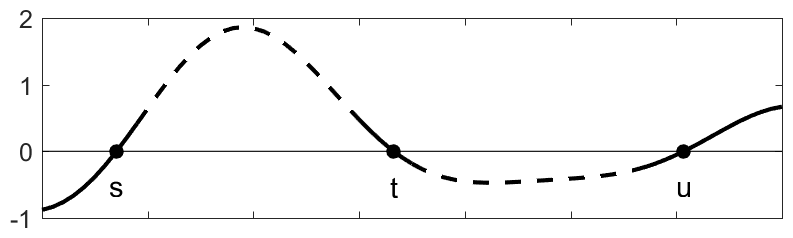}
}
\caption{ The ``triply conditioned'' process with successive crossing distances $T_1 = t-s, T_2 = u-t$, 
centred at $t$.}
\label{Fig:Ped}
\end{figure}

\section{Exact non-asymptotic crossing distributions}\label{s:exact}
\subsection{The Durbin-Rychlik formula}\label{ss:DRP}
The root of the exact formula for crossing interval distributions is the ``doubly conditioned'' 
process derived by \citet{Slepian1963}, which explicitly describes a process with a level upcrossing at  
$0$ and downcrossing at $t$; the conditioning is implicit in Rice's original paper \citep{Rice_b}.  
For our purpose, we illustrate in Figure~\ref{Fig:Ped} 
the ``triply conditioned'' process with a zero downcrossing at $t = 0$ with upcrossings at 
$s$ and $u$.   

We use the following  notational conventions from \citep{PodgorskiRM}: $X_{s,t,u} = (X(s), X(t), X(u))$, 
$\dot{X}_{s,t,u} = (\dot{X}(s), \dot{X}(t), \dot{X}(u))$, 
$\dot{X}_{s,t,u}^{+-+} = \dot{X}(s)^+ \dot{X}(t)^- \dot{X}(u)^+$, 
and $X_{s,t} = (u,v)$ means 
$X(s) = u, X(t) = v$. Moreover, $a \leq X_{s,t} \leq b$ 
means that for each $u \in (s,t): a \leq X(u) \leq b$, 
 while $\{a \leq X_{s,t} \leq b\}$ also stands for the indicator function of this set, 
i.e. the function equal to one whenever the condition between the brackets is 
satisfied and to zero otherwise. 
 
The triple crossing intensity for the configuration in Figure~\ref{Fig:Ped} is then equal to a truncated 
product moment in a conditional normal distribution:
\begin{align}
\nu^{+-+}(s,t,u)&= \int_{z_1 = 0}^\infty \int_{z_0 = -\infty}^0 \int_{z_2 = 0}^\infty 
z_1 z_0 z_2\, f_{\dot{X}_{s,t,u}, X_{s,t,u}}
(z_1, z_0, z_2, 0,0,0)\, \rd z_1\, \rd z_0 \, \rd z_2 \label{threederivatives} \\[0.4em]
&= \Ex [\dot{X}_{s,t,u}^{+-+} \mid X_{s,t,u} = (0,0,0)] \times
f_{X_{s,t,u}}(0,0,0). \label{triplecrossingintensity}
\end{align}
With $\nu =\int_0^\infty z f_{\dot{X}(0),X(0)}(z,0)\, \rd z$, the zero up/downcrossing 
intensity, $\nu^{+-+}(s,t,u)/\nu$ is the conditional intensity of upcrossings at $s, u$, given a 
downcrossing at $t$.
To get the distribution of successive crossing intervals one has to qualify the expectation in 
\eqref{triplecrossingintensity} 
by requiring that the process stays above $0$ in the entire left interval and below in the entire right interval in Figure~\ref{Fig:Ped}. 

The Rice series approximations achieves the qualifications by restricting the number of extra crossings in the interior 
of the intervals by higher order moments for the number of crossings. \citet{KanRobotti2017} give recursive  
formulas how to compute all truncated moments to a very high computational cost.   

The exact formula for the distribution of level crossing intervals in Gaussian processes was developed 
 by \citet{Durbin1985JAP} and \citet{Rychlik}, while \citet{PodgorskiRM} 
 extended it to successive intervals, without giving details. \citet[Thm.~7.1]{AbergRL} 
 later presented a complete proof for a very similar case; a short proof is given in Appendix~\ref{jointT1T2proof} 
 in the present work.
 
 Since we work with a stationary process, we can take $t=0$, and 
 consider the indicator functions to be included in the conditional expectation in \eqref{triplecrossingintensity}, 
 $ \{X_{s,0} > 0\}$  and  $\{X_{0,u} < 0\}$.   Obviously,  
 $\dot{X}_{s,0,u}^{+-+} = |\dot{X}_s \dot{X}_0 \dot{X}_u|$ when both conditions are satisfied.

The exact expression for the probability density of the length of two successive zero crossing 
intervals is, ($s = -t_1, u = t_2$), \citep[Eqn~10]{PodgorskiRM},
\begin{align}
&f_{T_1,T_2}(t_1, t_2)\nonumber \\ &= \nu^{-1} 
\Ex[|\dot{X}_{-t_1} \dot{X}_0 \dot{X}_{t_2}| \{X_{-t_1,0} > 0 > X_{0,t_2}\} \mid 
 X_{-t_1,0,t_2} = (0,0,0)]\, f_{X_{-t_1,0,t_2}}(0,0,0). \label{jointpdf}
\end{align}
In the above, the conditional expectation is taken over infinite dimensional set of 
variables due to the uncountable  number of times instants involved 
in $\{X_{-t_1,0} > 0 > X_{0,t_2}\}$.

\subsection{Evaluation of the expectation involving an uncountable number of instants}
Like most, so called, ``explicit solutions'' to mathematical problems, 
the expectation in \eqref{jointpdf} has to be evaluated 
numerically.\footnote{The sad truth is that there are no known closed forms of the 
Gamma function for irrational values. (StackExchange)}
The degree of complexity is the same as computing the distribution of the maximum of a smooth 
non-stationary Gaussian process $X(t)$, $\Prob (\max_{[0,T]} X(t) \leq x)$,  
for a finite interval with length $T$; see  
 \citep{GenzBretz2009} for available efficient software, and \citep{AzaisGenz2013} for an 
 analysis of numerical accuracy.  

For our problem, we observe that the variables in the braced indicator in \eqref{jointpdf} 
are non-stationary Gaussian, given the condition $ X_{-t_1,0,t_2} = (0,0,0)$. 
The derivatives at the crossing points are not Gaussian but their joint 
density is proportional to the integrand in \eqref{threederivatives}. If we incorporate the derivatives 
in the conditioning, the indicator variables are still non-stationary Gaussian.    

Due to the strong local dependence for smooth Gaussian processes,  
the natural way to compute the expectation in \eqref{jointpdf} is to 
replace the ``infinite-dimensional'' indicator  $\{X_{-t_1,0} > 0 > X_{0,t_2}\}$
by a finite-dimensional one. To obtain sufficient accuracy one may have to take a dense grid which can 
result in an almost singular covariance matrix for the multivariate Gaussian distribution. 
\citet{Brodtkorb2004,Brodtkorb2006} discussed several strategies to evaluate nearly singular multinormal 
expectations. He improved the algorithms proposed by \citet{Genz1992} and 
\citet{GenzKwong2000} when the correlation is strong and the number of variables is very large, and, 
most important, he increased the computing speed and improved memory requirement by utilizing ideas 
from \citep{Rychlik1987AAP,Rychlik1992,PodgorskiRM}. 
\citet{Brodtkorb2006} also made extensive studies of the accuracy of the numerical algorithms 
for a number of realistic applications. The result, the {\sc Matlab} routine {\sf RIND} is included in 
the free package {\sf WAFO}, \citep{gitwafo}, and in the {\sf MAGP} package by \citet{Mercadier}..

\subsection{About {\sf RIND}}\label{s:aboutrind}

The {\sf RIND} routine is  designed to accurately approximate functionals like
\begin{align}
\Ex[|\dot{X}_{-t_1} \dot{X}_0 \dot{X}_{t_2}| \{X_{-t_1,0} > 0 > X_{0,t_2}\} \mid 
 X_{-t_1,0,t_2} = (0,0,0)] f_{X_{-t_1,0,t_2}}(0,0,0), \label{Rindexpectation}
\end{align}
when $ \{X_{-t_1,0} > 0 > X_{0,t_2}\}$ is replaced by a discrete time restriction, 
 $ \{X_{S_n'} > 0 > X_{S_m''}\}$, on equidistant 
points $s_k', s_k''$ in the two intervals: 
\begin{equation}
\begin{split}
S_n' &= (-t_1 < s_n' < s_{n-1}' < \ldots < s_1' < 0),\\
S_m'' &= (0 < s_1'' < s_2'' < \ldots < s_m'' < t_2).
\end{split}
\label{discreteS}
\end{equation}
Thus, the problem falls in the category of general multinormal expectations: 
\begin{align}
F(\bol{a}, \bol{b}; \bol{\Sigma}) &= \frac{(2N)^{-N/2}}{\sqrt{|\bol{\Sigma}}|} \int_{a_1}^{b_1} 
\cdots \int_{a_N}^{b_N} g(\bol{x}) \exp \left\{\frac{-\bol{x}^T \bol{\Sigma}^{-1} \bol{x}}{2} \right\}
\, \rd x_N \cdots \rd x_1 .
\end{align}
The paper by \citet{Genz1992} is the main reference to modern techniques for evaluation of the
integral. Its focus is on the case $g(\bol{x}) = 1$, but it gives hints on how to handle a general 
$g$-function. The work by \citet{Brodtkorb2006} is focused on the special structure of 
\eqref{Rindexpectation} in order to increase efficiency, while using the basic ideas from \citep{Genz1992} 
and subsequent work. The additional elements include use of the regression approximation 
 by \citet{Rychlik1987AAP}, removal of redundant integrals, and Cholesky matrix truncation.

The result, the {\sf RIND} routine, offers several alternative methods for the computation of the 
integral, including combinations of the following codes by Brodtkorb (2000 and 2004); more details 
of the alternatives can be found in the {\sf Fortran} source file {\sf intmodule.f} in {\sf WAFO}. 
\begin{description}\setlength\itemsep{0mm}
\item[SADAPT:] A generalization by  Brodtkorb (2000) of the routine {\sf SADMVN} in \citep{Genz1992} 
to make it work not just for the multivariate normal integral. 
\item[KRBVRC:] An update by Brodtkorb (2000) of the module {\sf KRBVRCMOD} by Genz (1998).
\item[KROBOV:] An update by Brodtkorb (2000) of the module {\sf KROBOVMOD}  by Genz (1998).
\item[RCRUDE:] An update by Brodtkorb (2000) of the module {\sf RCRUDEMOD} by Genz (1998), improving randomized integration.
\item[SOBNIED:] A routine by Brodkorb (2004), improving 
{\sf KRBVRC} by using random selection of points as in \citep{HongHickernell2003}.
\item[DKBVRC:] An update by Brodtkorb (2004) of a routine with the same name by Genz (2003). 
\end{description} 
 
\subsection{Calling {\sf RIND}}
 The {\sf RIND} function is a {\sc Matlab} interface to a set of algorithms, originally written 
 in {\sf Fortran} and {\sf C++}, that execute one, or a combination, of the options, {\sf SADAPT -- DKBVRC}.   
The function takes as input means and covariances of three groups of normal variables. 
One groups consists of variables to condition on, {\sf Xc = xc}, like $X_{-t_1,0,t_2} = (0,0,0)$ in 
\eqref{Rindexpectation}. The second group, {\sf Xd}, are the derivatives at crossings, like 
$\dot{X}_{-t_1}, \dot{X}_0, \dot{X}_{t_2}$, and the third group {\sf Xt} contains the variables 
that have to satisfy an interval condition, {\sf xlo} $\sfleq$ {\sf Xt} $\sfleq$ {\sf xup}. Such constraints 
may also be imposed on the derivatives. 

The call to {\sf RIND} has the following input/output structure, extracted from {\sf help RIND}: 

\bigskip
\begin{minipage}{0.8\textwidth}\small
\begin{verbatim}
[E,err,terr,exTime,options] = rind(S,m,Blo,Bup,indI,xc,Nt,options);
         E = expectation/density, according to (4)
       err = estimated sampling error
      terr = estimated truncation error.
    exTime = execution time
         S = Covariance matrix of X=[Xt;Xd;Xc]
         m = the expectation of X=[Xt;Xd;Xc] 
   Blo,Bup = Lower and upper barriers used to compute the integration bounds
      indI = indices to the different barriers in the indicator function
        xc = values to condition on 
        Nt = size of Xt
   options = options structure or named parameters with corresponding values
\end{verbatim}
\end{minipage}
\bigskip

The options structure is used to select integration method, set error tolerances, alternatively set 
speed, set seed for Monte Carlo-integration, and many other parameters.  
As an example of how the {\sf speed} option affects the result we can take the low-frequency white 
noise process, illustrated in Figure~\ref{Fig:IIA}, {\sf WN}, when integrated by {\sf SOBNIED} at 
slowest and fastest speed, {\sf speed = 1}, and {\sf 9}, respectively. The joint pdf of two 
successive zero crossing intervals was calculated at a grid of {\sf 120 x 120 = 14 400} points. 
Execution time was 120 seconds with the slowest option and  45 seconds 
with the fastest, with an increase in truncation error by a factor $10^{4}$  and in  
sampling error by a factor $4$ to $10^{2}$. The level curves were virtually identical 
for levels enclosing up to $90$\% of the probability, with only small deviations for 
the more extremes curves. The plot in Figure~\ref{Fig:IIA} was produced in 81 seconds with the 
{\sf KROBOV} algorithm, a generally slower method, with {\sf speed} set to 5.

\begin{remark}[{\sf RIND} and persistence]~\label{RINDandPers} 
The {\sf RIND} function can be used to compute the 
persistence for any finite interval $[0, T]$. One just has to let the groups 
{\sf Xc} and {\sf Xd} of conditioning variables and 
end point derivatives be empty. Setting the lower and upper bounds 
to $0$ and $\infty$, respectively, the algorithm 
will give the probability that all the {\sf Xt}-variables stays above $0$. 
The covariance matrix {\sf S} has to be set to the covariance matrix
 of a sufficiently dense subset of $X(t)$-variables, $X(s_k); 0 < s_k < T$. 
\end{remark}

\section{{\sf RIND} and successive level crossing intervals}
In this section we give an account of numerical evaluation of crossing distribution properties for 
Gaussian processes. A complete code using the {\sf WAFO} toolbox for some of the evaluations is provided in 
Appendix~\ref{code} to illustrate convenience and simplicity of using the implemented integration routines. 
It was developed for use on joint crossing intervals in \citep{Lin2019}.
 
\subsection{Selection of spectra and covariance functions}
We study throughout a stationary Gaussian process with mean zero and 
covariance function $r(t) = \Cov (X(s), X(s+t))$. All processes are smooth 
with continuously differentiable sample paths and finite number of level crossings 
in bounded intervals. We refer to \citep{Lindgren2013} for general facts 
on stationary Gaussian processes, and for differentiability in particular.

We illustrate the {\sf RIND} method on processes with covariances and spectra 
of many different types, studied in the literature:  
\citep{LonguetHiggins1962,Lindgren1972AAP,Sire2008}, 
\citep{AzaisW,BrayMajumdarSchehrAiP2013,WilsonHopcraft2017}, all listed in Table~\ref{Tab:1}. 

\begin{table}
\caption{Un-normalised spectra and covariance functions.}
\vspace{2mm}
\label{Tab:1} 
\footnotesize
\centerline{
\begin{tabular}{lll}\toprule
\sc Spectrum 
& \sc Covariance function & \sc Source \\ 
\midrule[0.9pt]
\multicolumn{3}{c}{Type: Rational spectrum}  \\
\midrule 
$(1 + \omega^2)^{-3}$ & $e^{-|t|}(1 + |t| + t^2/3)$ & {\sf LH1} \\
$(1 + \omega^2)^{-4}$ & $e^{-|t|}(1 + |t| + 6t^2/15+|t|^3/15)$ &  {\sf LH2} \\
$(1 + \omega^2)^{-5}$ & $e^{-|t|}(1+|t|+3t^2/7 + 2|t|^3/21 + t^4/105)$ & {\sf LH3} \\
${\omega^4}{(1+\omega^2)^{-5}}$ & $e^{-|t|}(1 + |t| - t^2/3 - 2|t|^3/3 + t^4/9)$ & {\sf LH4}  \\
$(1 + \omega^2)^{-2}$ & $e^{-|t|}(1 + |t|)$ & {\sf LH5} \\
${\omega^2}{(1+\omega^2)^{-4}}$ & $e^{-|t|}(1 + |t| - t^2/3)$ & {\sf LH6} \\
${\omega^4}{(1+\omega^2)^{-4}}$ & $e^{-|t|}(1+|t|-2t^2+|t|^3/3)$ &  {\sf LH7} \\
\midrule[1pt]
\multicolumn{3}{c}{Type: Shifted Gaussian} \\ 
\midrule
$\cosh (k\omega)\exp (-\omega^2/2)$ & $\cos (k t)  \exp(-t^2/2), k = 0,1,2,...$ & {\sf WHk}\\
\midrule[1pt]
\multicolumn{3}{c}{Type: Noise and sea waves}\\ 
\midrule
$\bone_{[-1, 1]}$ & $\sin(t)/t$ & {\sf WN} \\
$(1+\omega^{14})^{-1}$ & NA & {\sf BS} \\
Jonswap & NA & {\sf J} \\
\midrule[1pt]
\multicolumn{3}{c}{Type: Diffusion} \\ 
\midrule
$\displaystyle \operatorname{sech}(\pi \omega)
=\frac{1}{\cosh(\pi \omega)},~d=2$&$\displaystyle{\operatorname{sech}^{d/2}(t/2)
=\frac{1}{\cosh^{d/2}\left(t/2\right)},~d\in \mathbb N}$&{\sf  BMSd}\\
\bottomrule 
\end{tabular}
}
\end{table}

Most of the studied processes are, what was called by \citet{LonguetHiggins1962}, 
the ``regular type'', with covariance function admitting the expansion 
\begin{align} 
r(t) &= 1 - \lambda_2 t^2/2 + \lambda_4 t^4/4! + o(t^{4+\epsilon}), t \to 0. \label{regular}
\end{align}
In the regular case, the process is twice continuously differentiable and the zero crossing interval 
pdf is of order $f_T(t) \thicksim c t$ for small $t$. The zeros do not cluster.

Processes of the ``irregular type'' have covariance functions expanded as 
\begin{align}
r(t) &= 1 - \lambda_2 t^2/2 + C |t|^3/3! +  o(|t|^3), C \neq 0, \label{irregular}
\end{align}
and the crossing density has a universal non-zero limit at the origin, 
\begin{align}
f_T(t) \to K \alpha, \text{ where } \alpha = C/6\lambda_2. \label{irregularlimit}
\end{align}
\citet[Eqn~47]{LonguetHiggins1963} gives an estimate of the constant $K \approx 1.15597$. 
 
The name convention for the spectra and covariances in Table~\ref{Tab:1} is as follows. 
{\sf LHk} and {\sf WHk} hint at \citep{LonguetHiggins1962} and \citep{WilsonHopcraft2017}. The 
diffusion spectra {\sf BMSd} were used in \citep{BrayMajumdarSchehrAiP2013}. 
{\sf WN} is low-frequency white noise with a Butterworth approximation {\sf BS}, and the 
Jonswap spectrum {\sf J} is an example of an ocean wave spectrum. 
Note that the spectra and covariance 
functions are listed in un-normalised form. In the examples they are normalised to 
$\lambda_0 = \lambda_2 = 1$ and average zero crossing interval equal to $\pi$.  
All spectra are of the regular type 
except {\sf LH5--LH7}, which are irregular. 

The regularity parameter is $\alpha = \lambda_2/\sqrt{\lambda_0 \lambda_4}$; $1/\alpha$ is equal to 
the mean number of local extremes per mean level crossing. In the numerical examples we will relate 
the $\alpha$-parameter to the correlation coefficient between successive crossing intervals, and to the 
deviation between the true joint pdf and the pdf under the Independent Interval Assumption, 
obtained by multiplication of the marginal pdf:s. 
The Kullback-Leibler distance, 
$$KL = - \iint f_{T_1,T_2}(t_1,t_2) 
\log (f_{T_1}(t_1)f_{T_2}(t_2)/f_{T_1,T_2}(t_1,t_2))\, \rd t_1 \rd t_2 ,$$ 
is used as a measure of the deviation between full dependence and independence. 
Note that we use the {\sl true marginal pdf\,} for a Gaussian process when we construct the 
``independence'' pdf, and not the one obtained by the IIA technique. 
One of the reasons of the difficulty in using the latter is, as explained in 
Appendix~\ref{Val-IIA}, that the IIA approach typically does not produce a valid probability distribution.

\subsection{Computations of the joint distribution of crossing intervals}
The joint pdf $f_{T_1,T_2}(t_1, t_2)$ for successive mean level crossing intervals are 
calculated by {\sf RIND} as described in Section~\ref{s:aboutrind} via the user 
interface {\sf cov2ttpdf}. It takes as input the covariance function in the form of a {\sc Matlab} 
symbolic function or the spectrum in the form of a {\sf WAFO} spectrum structure. 
Both types of arguments are automatically normalized to $\lambda_0 = \lambda_2 = 1$. 
For irregular type processes, the form of the density for small $t_1, t_2$ is not directly 
resolved by the integration in {\sf RIND} but the limit \eqref{irregularlimit} 
can be used to extend the integrated values to $t_1=0$ and $t_2=0$. 
The plot for {\sf LH7} in Figure~\ref{Fig:nonIIA} is constructed in that way. 
The marginal pdf of the crossing intervals $T_1$ and $T_2$ are obtained by integrating 
the bivariate pdf, and a joint pdf under independence is obtained by multiplying the two.

\begin{figure}
\centerline{
\includegraphics[height=70mm]{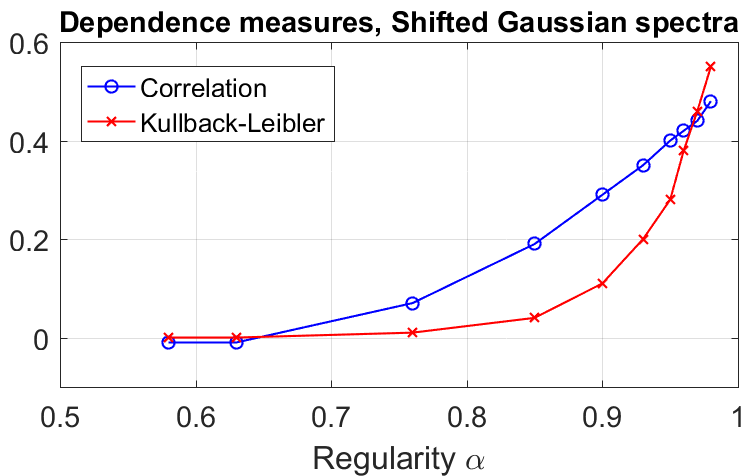}
}
\caption{ Dependence measures for shifted Gaussian spectra.}
\label{Fig:Measures}
\end{figure}

We have compared the exact {\sf RIND} pdf with the IIA-based pdf for all the spectra in 
Table~\ref{Tab:1} and we present some typical illustrations. We also give some 
numerical dependence measures: in  Figure~\ref{Fig:Measures} for shifted Gaussian spectra and in Table~\ref{Tab:othermeasures} for other  examples of spectra.

For each of the spectra, we computed the correlation coefficient between successive half periods, 
and the Kullback-Leibler distance between the exact pdf and the IIA pdf. Figure~\ref{Fig:Measures} shows the 
smooth relation between the regularity $\alpha$ and the dependence for the shifted Gaussian spectra, 
{\sf WH0-WH9}. The theoretical correlations presented in Figure~\ref{Fig:Measures} 
agree with those illustrated in \citep[Fig.~8]{WilsonHopcraft2017}. 
The spectra in the other group are more diverse and do not exhibit any systematic 
relation with the regularity measure, (when it exists),  as seen in Table~\ref{Tab:othermeasures}. 

\begin{table}[h!]
{\tabulinesep=1.15mm
\caption{Dependence measures for rational, noise, and wave spectrum}
\label{Tab:othermeasures}
\vspace{2mm}
\centerline{\setlength\tabcolsep{3.8pt}
\begin{tabular}{lcccccccccc} \toprule
model  & {\sf LH1} & {\sf LH2} & {\sf LH3} & {\sf LH4} & {\sf LH5} & {\sf LH6} & {\sf LH7} %
& {\sf WN} & {\sf BS} & {\sf J}  \\ \midrule
$\alpha$ 				& 0.33 & 0.45 & 0.49 & 0.49 & NA & 0.50  & NA & 0.74 & 0.72 & 0.70 \\
$\Corr (T_1,T_2)$ 	& -0.02 & -0.02 & -0.02 & 0.25 & 0.04 & 0.19 & 0.23 & -0.01 & -0.02 & 0.44 \\
{\sf KL}              		& 0.00 & 0.00 & 0.00 & 0.03 & 0.01 & 0.02 & 0.05 & 0.04 & 0.01 & 0.14 \\ \bottomrule
\end{tabular}
}
}
\end{table}

We illustrate graphically the results for a selection of different spectra. 
Each plot in Figures~\ref{Fig:SH1210joint}-\ref{Fig:nonIIA} shows three sets of bivariate distributions, 
illustrated by level curves enclosing 10, 30, \ldots, 99.9~\% of the distributions: 
red curves for the exact pdf, blue for the synthetic independent pdf, and black dashed 
for simulated data with about 2.6~million crest-trough interval pairs.  
The smooth appearance of the blue curves is due to the integration. 
We show four examples with almost independent half periods, and four examples 
with very clear and diversified dependence.

\subsubsection{Gaussian diffusion {\sf BMSd} in dimension $d$}\label{sss:diffusion}
We start with a class of stationary Gaussian processes where successive zero-crossing intervals are 
``almost'' independent, representing diffusion in $d$ dimensions, where $d = 1, 2, 3$ 
correspond to physically feasible experiments, \cite[Sec.~9]{BrayMajumdarSchehrAiP2013}. 
The covariance function is $r(t) = \sech^{d/2}(t/2)$. We used the {\sf RIND} function with 
the {\sf SOBNIED} integration routine to compute the joint pdf \eqref{jointpdf} of two successive 
zero-crossing intervals. The speed parameter was set to $1$, which gives the most accurate results. 
With a resolution in \eqref{discreteS} of $\Delta s = 0.2$ and a total of $125$ points for both crest 
and trough periods the computation time was 260 seconds for each diagram. 

\begin{figure}
\centerline{
\includegraphics[width=45mm]{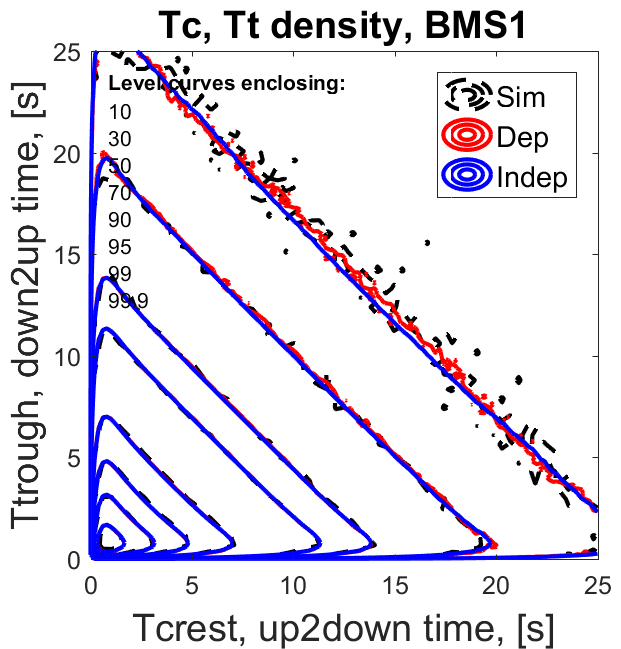}
\includegraphics[width=45mm]{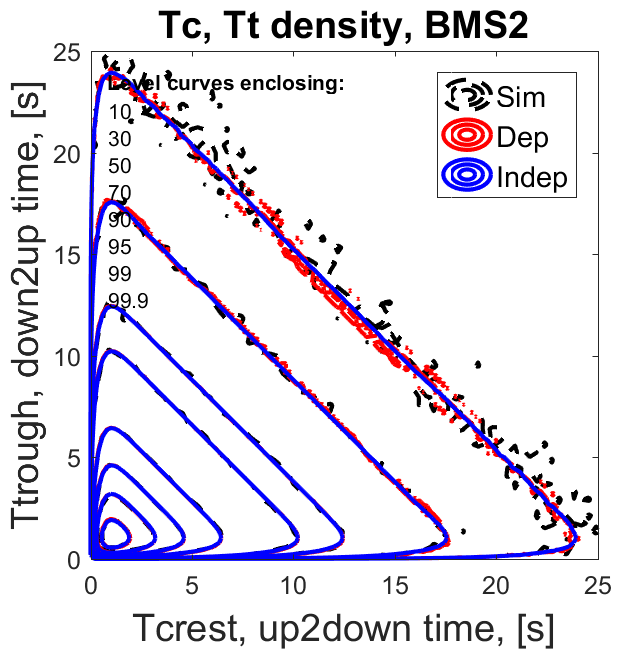}
\includegraphics[width=45mm]{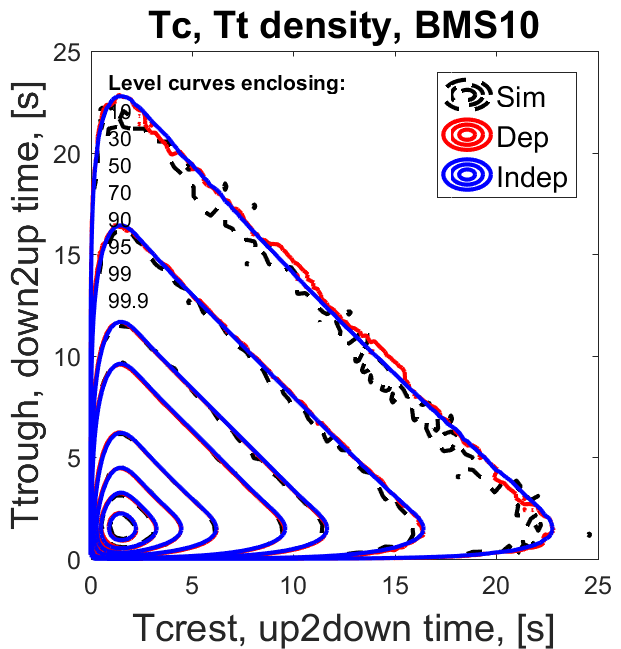}
}
\caption{Joint pdf for crest-trough periods for Gaussian processes with with 1D, 2D, 10D  diffusion 
type covariance {\sf BMSd}.}
\label{Fig:SH1210joint}
\end{figure}

\begin{figure}
\centerline{
\includegraphics[width=40mm]{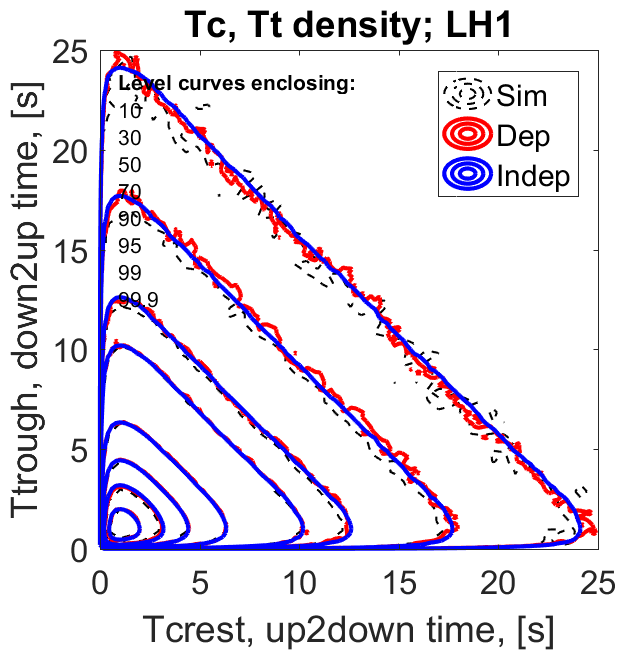}
\includegraphics[width=40mm]{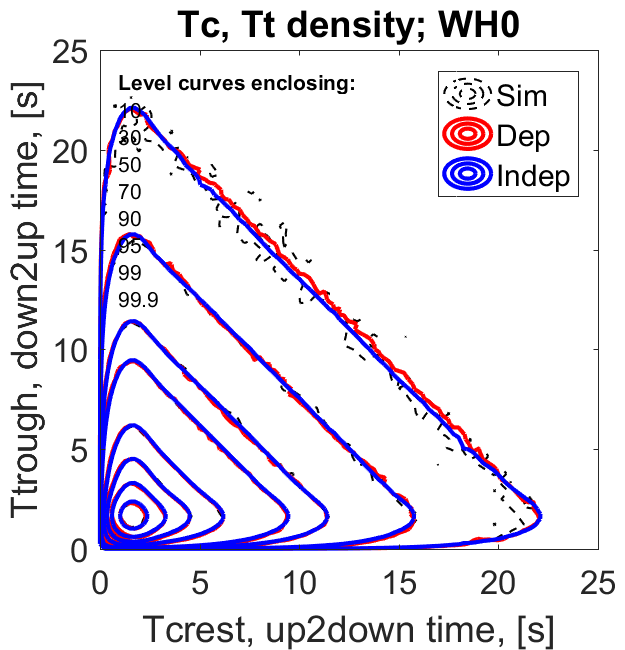}
\includegraphics[width=40mm]{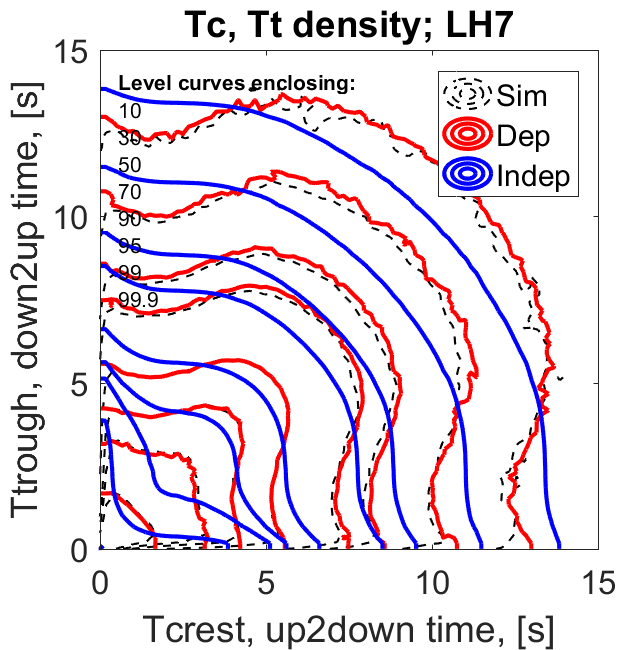}
\includegraphics[width=40mm]{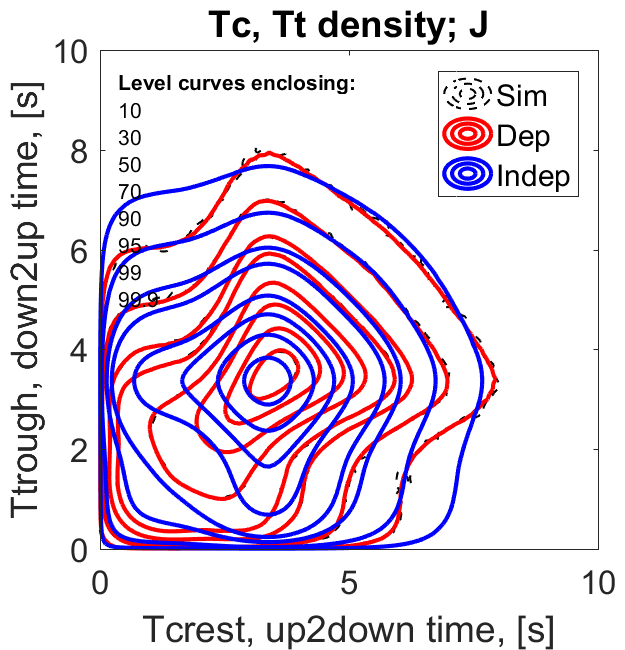}}
\vspace{3mm}
\centerline{
\includegraphics[width=40mm]{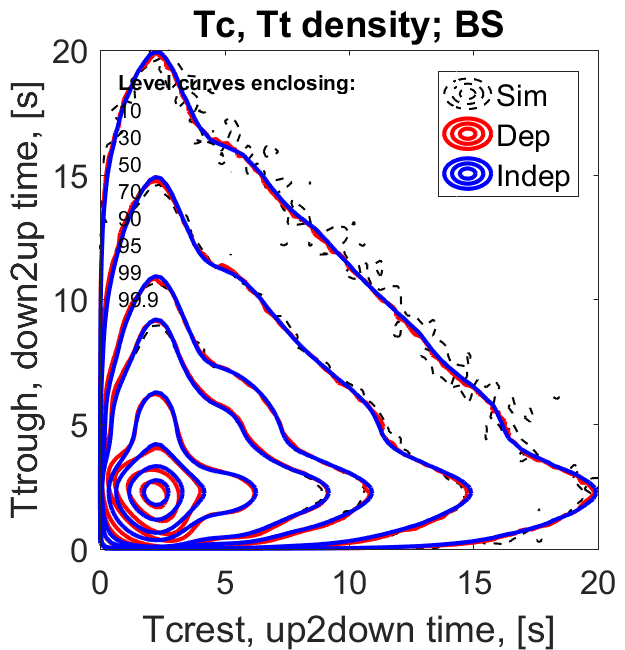}
\includegraphics[width=40mm]{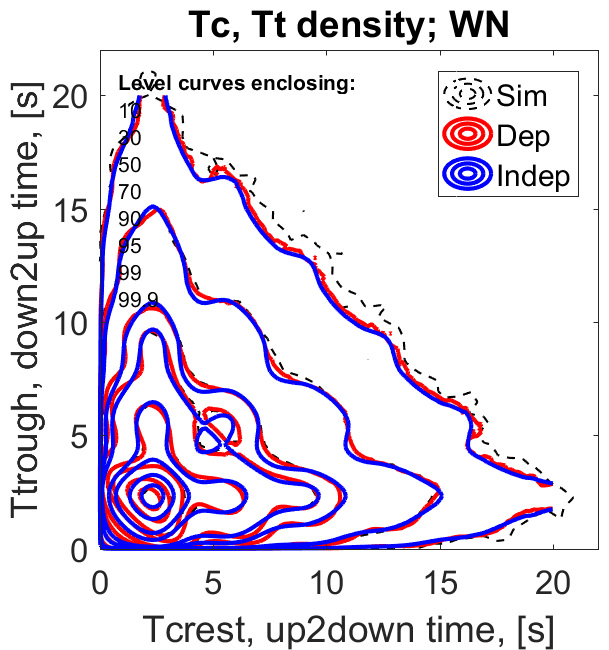}
\includegraphics[width=40mm]{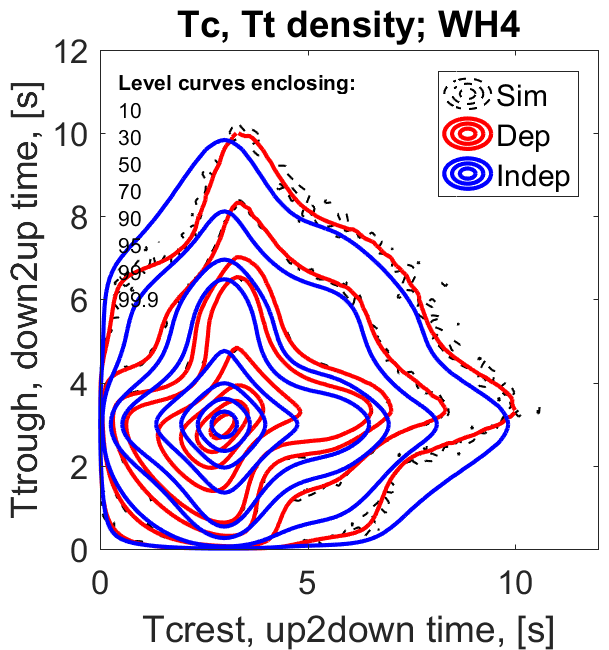}
\includegraphics[width=40mm]{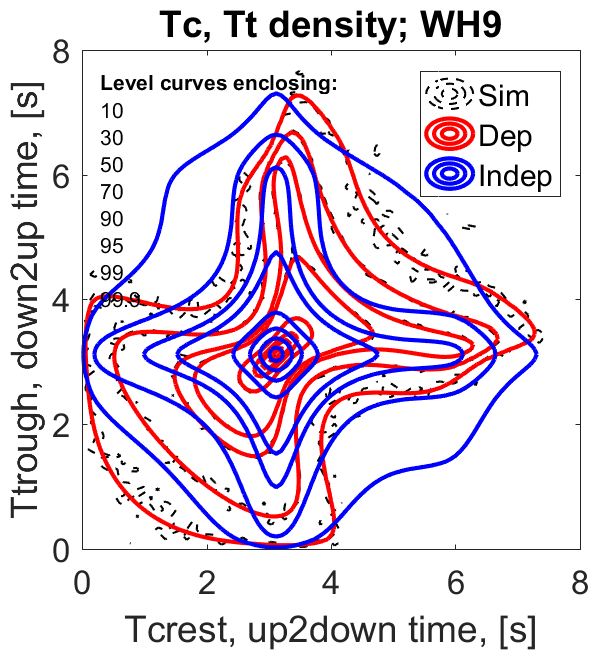}
}
\caption{Joint pdf for spectra.  The first two columns (left) almost independent half periods: {\sf LH1, WH0, BS, WN} with correlation coefficient {\sf -0.02, -0.02, -0.02, -0.01} and  Kullback-Leibler distance {\sf KL = 0.00, 0.00, 0.01, 0.04}, respectively.  
The remaining two columns (right) spectra with clear dependence: : {\sf LH7,  J, WH4, WH9} with correlation coefficient 
{\sf 0.13, 0.44, 0.30, 0.48} and Kullback-Leibler distance {\sf 0.05, 0.14, 0.11, 0.55}, 
respectively.
Level curves enclose 
10, 30, \ldots, 99.9~\% of the distributions. Red solid curves: pdf computed by {\sf RIND}; 
Blue solid curves: pdf under IIA assumption, with marginal pdf from {\sf RIND}; 
Black dashed curves: simulated pdf from about 2.6~million pairs of crest-trough intervals.}
\label{Fig:IIA}
\end{figure}

The results are shown in Figure~\ref{Fig:SH1210joint}. 
For all quantile levels, except for the most extreme, the true pdf and the one obtained by 
the independent approximation are almost identical. The pdf:s are very favourably compared to 
empirical histograms, based on about 2.6 million crest-trough period pairs each. We return to this 
example in Section~\ref{ss:BMS}.

\subsubsection{Other spectra with near independent half periods}\label{sss:otherindep}
Figure~\ref{Fig:IIA} shows joint pdf for four different type spectra where successive intervals 
are near uncorrelated and almost pairwise independent, as measured by the 
Kullback-Leibler distance.   There is clear visual 
agreement between the exact pdf (red) and the pdf with independent margins. 
The Gaussian spectrum {\sf WH0} and the approximating {\sf LH1} are centered at zero frequency 
and follow the IIA model almost perfectly.  
However,  the white noise spectrum {\sf WN} deviates considerable from how one 
normally  envisages independent variables. The Butterworth spectrum {\sf BS} approximates the 
{\sf WN} spectrum and the pdf is near to independence.
Note again that the 
marginal distributions are computed from the exact pdf  by integration and not by the renewal  
argument as in the IIA approach.  

\subsubsection{Spectra with strongly dependent half periods}\label{sss:dependent} 
Figure~\ref{Fig:nonIIA} shows examples with clear or even strong dependence. The spectrum 
{\sf LH7} is of the irregular type \eqref{irregular} with a small tendency of having three zeros close to 
each other, as shown by the red exact level curve near the origin. The three other spectra, 
the {\sf Jonswap, J} and the shifted Gaussians, {\sf WH4, WH9},   are regular 
with large regularity parameter $\alpha$, large correlation coefficients, and Kullback-Leibler 
difference larger than $0.1$. As is obvious from the figures the dependence can take many 
different shapes, which makes it difficult to catch it in a simple parametric form.

\subsection{Gaussian diffusion {\sf BMSd} and the persistence exponent}\label{ss:BMS}
The diffusion type covariance appears as a model for the diffusive time development of 
a Gaussian random field, initiated as white noise at time $t=0$. After a transformation to 
logarithmic time, $T = \log t$, the field is a stationary (homogeneous) Gaussian field with 
the $\sech^{d/2}(t/2)$ covariance function.    
Persistence for diffusion systems has been studied in physics since the early 1990s.  

\cite{MajumdarSBC} used the Independent Interval Assumption and \eqref{ThetaIIA} to 
find $\theta_{IIA}$ for different dimensions, and compared with simulations, while 
\cite{WongMairWalsworthCory2001} presented experimental evidence for $d=1$. 
\cite{NewmanLoinaz2001} designed an efficient simulation procedure to estimate the persistence probability for arbitrary dimension and suggest corresponding persistence exponents based on 
the simulations. 
\cite{PoplavskyiSchehr2018} gives a definite answer for dimension 
$d=2$, namely $\theta(2) = 3/16 = 0.1875$, a value that agrees with the  
\cite{NewmanLoinaz2001} simulation. At present, no exact values are known for other dimensions. 

We can now compare these 
results with values computed by means of the {\sf RIND}, which integrates the 
multidimensional normal density to give the tail probability $Q_T$ \eqref{QT} 
for arbitrary interval length; Remark~\ref{RINDandPers}.  

Our first concern is the shape of the distribution. 
In Section~\ref{sss:diffusion} we argued that the independence of successive intervals appears to be approximately satisfied for the diffusion spectra, Figure~\ref{Fig:SH1210joint}. The marginal distributions 
are rather close to exponential, even if not exactly so. What effect the deviation has on the 
persistence approximation $\theta_{IIA}$ from \eqref{ThetaIIA} is unclear. 
However, as shown in Appendix~\ref{ddt}, Table~\ref{perco}, and with agreement with earlier results in the literature, the resulting approximation of the persistency exponent is in the vicinity of 0.1863 which underestimated the true value 0.1875.  Furthermore, 
the structure of the theoretical limiting behaviour \eqref{exponentialindex} makes it difficult 
to estimate any persistence exponent by simulation or numerical computation.

\begin{table}
{\tabulinesep=1.15mm
\caption{Numerically calculated persistence exponents $\theta(d)$ for diffusion. 
$^*$For the values in parenthesis, see text.}
\label{Tab:diffusionexponents}
\vspace{3mm}
\centerline{\setlength\tabcolsep{3.8pt}
\begin{tabular}{ccc|ccc} \toprule
d & NL & RIND   & d & NL & RIND\\ \midrule
1 & 0.1205        & 0.1206 (0.1203)$^*$ & 10 & 0.4587 & 0.4589\\
2 & {\bf 0.1875} & 0.1874 (0.1875)$^*$ & 20 & 0.6556 & 0.6561\\
3 & 0.2382        & 0.2382 & 30 & 0.8053 &  0.8063\\
4 & 0.2806        & 0.2805 & 40 & 0.9232 & 0.9327\\
5 & 0.3173        & 0.3171 & 50 & 1.0415 & 1.0439\\
\bottomrule
\end{tabular}
}
}
\end{table}

The persistence estimates in \citep[Tab.~1]{NewmanLoinaz2001} come from a large and 
well controlled simulation experiment of over $10^8$ realizations of first crossing events. We 
have designed the following numeric procedure to estimate the persistence exponent. 
\begin{enumerate}\setlength\itemsep{-1mm}
\item Fix a maximum $T_m$ and a grid $t_k = k\Delta$, $t_m \Delta = T_m$; 
this resembles the logarithmic grid in the \cite{NewmanLoinaz2001} simulations. $T_m$ 
should be chosen so large that any numeric or stochastic uncertainty in the tail is revealed.  
\item Compute, with {\sf RIND}, $Q_T$ for $T = t_k, k=1,\ldots, m$;  \cite{NewmanLoinaz2001} 
estimate $Q_T$ by simulation.
\item Make repeated independent runs to compute $Q_T$ and take the average $\overline Q_T$; 
{\sf RIND} uses Monte Carlo 
integration for extreme cases, and averaging will reduce the stochastic uncertainty. 
\item Make a regression of $-\log (2*\overline Q_T)$ against $T$; 
this can be done globally, over the entire 
interval $[0, T_m]$ or locally to allow for slow asymptotics in \eqref{exponentialindex}. 
The persistence index is then estimated as the slope of the fit. 
We used  quadratic polynomial fit to identify any typical trend in the exponential 
$e^{(\theta + A(T))T}$. 
\end{enumerate}  

We used this scheme to estimate the persistence index for the 10 selected 
dimensions in \citep{NewmanLoinaz2001}. 
We set $\Delta = 0.05$ and $T_m = 15$ for $d \geq 4$, which corresponds to the time span used in 
that paper and $\Delta = 0.1$ and $T_m = 30, 30, 20$ for $d = 1, 2, 3$. 
We computed $\overline Q_T$ as the average of 400 independent runs of 
{\sf RIND} with the {\sf SOBNIED} method and highest precision. The quadratic fit indicated 
that the local slope increased slowly with $T$. 
Table~\ref{Tab:diffusionexponents} shows the local $\theta_T$ in the middle of the interval, 
$T = T_m/2$. 

For dimensions $d=1, 2$ we give two values in parenthesis to take account of the asymptotic character 
of the persistence exponent. To accomplish this, we increased $T_m$ to  $40, 35$, respectively 
and estimated the local rate of decay for large $T$.  For $d=1$, the value $\theta = 0.1206$ is the best 
estimates over the interval $(0,30)$, while the value $0.1203$ is the stable value for large $T$. 
The {\sf NL}-value $0.1205$ is probably too high, based on a too short time span. The {\sf RIND}-value 
$0.1203$ agrees with the value computed by IIA. 
For $d=2$ the stable value is $0.1875$, equal to the theoretical value $3/16$, \cite{PoplavskyiSchehr2018}.

\begin{figure}
\centerline{
\includegraphics[height=75mm]{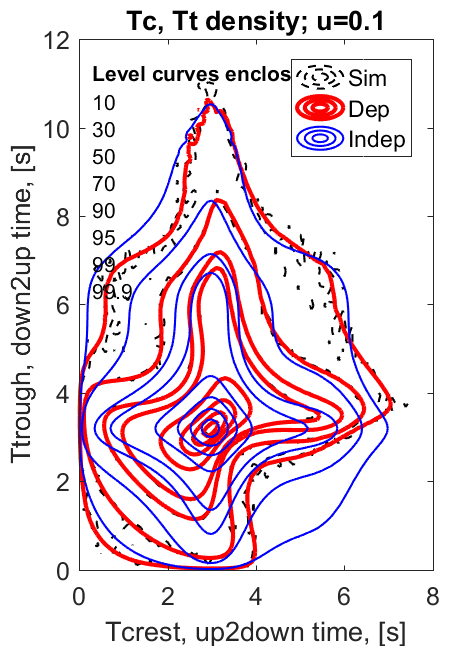} \hspace{3mm}
\includegraphics[height=75mm]{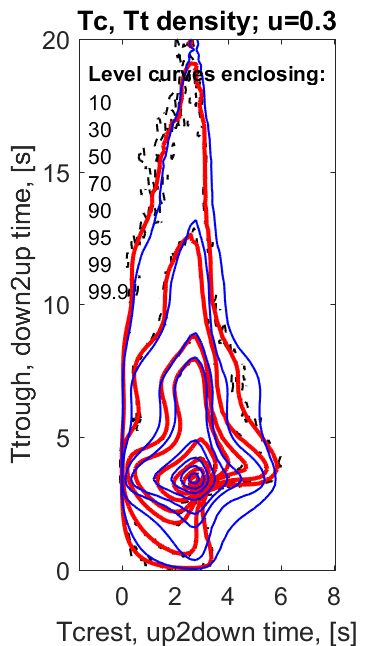} \hspace{3mm}
\includegraphics[height=75mm]{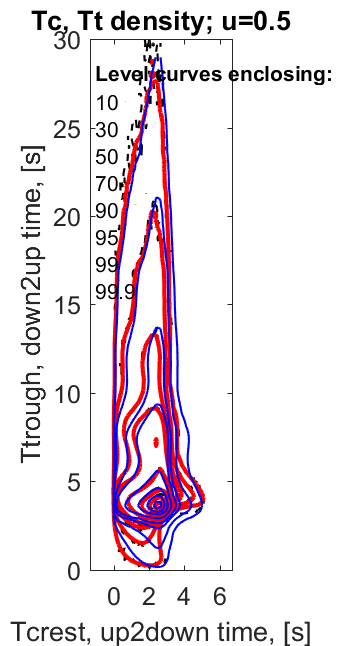}
}
\caption{Joint pdf of interval above and below levels $u = 0.1, 0.3, 0.5$ for shifted Gaussian 
spectrum {\sf WH6}.} 
\label{Fig:WH6u}
\end{figure}

\subsection{Crossings of non-zero levels}
The {\sf RIND} function is not limited to joint crest-trough period distribution but can be used 
for crossings of any level. We illustrate this on the Shifted Gaussian spectrum model 
{\sf WH6} for levels $u = 0.1, 0.3, 0.5$. Figure~\ref{Fig:WH6u} shows the dependence, 
and also the good agreement between the {\sf RIND} results and simulated pdf:s, 
based on more than 2.2~million pairs of excursions. 

\begin{figure}
\centerline{
\includegraphics[height=80mm]{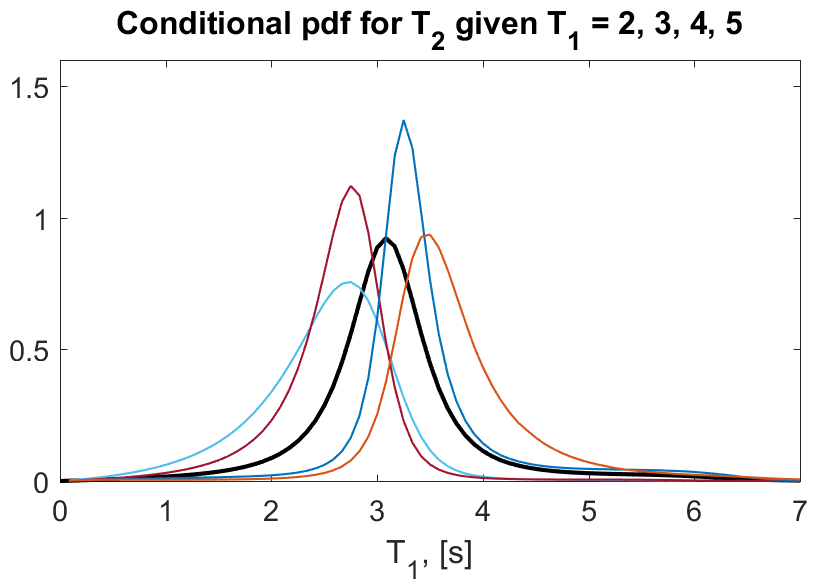}
}
\caption{Crossing and conditional crossing densities for WH6. 
Black curve: marginal pdf for zero crossing interval; colour curves: 
conditional pdf for next interval given the first one, $T_2 \mid T_1 = 2, 3, 4, 5$.}
\label{Fig:ConditionalWH6} 
\end{figure}

\subsection{The Markov approximation}\label{ss:Markov}
The Markov approximation for zero crossing intervals dates back to \citet{McFadden1958} 
and \citet{Rainal1962,Rainal1963}, who also tested the model by means of the sequence of 
correlations. 
Having the exact joint distribution of successive zero crossing intervals a natural next step is 
to test the Markov chain dependence of the full sequence of crossing intervals. If 
$f_{T_1,T_2}(t_1,t_2)$ is the joint density of two successive crossing intervals, one can construct 
a Markov transition kernel as $k(t_2 \mid t_1) = f_{T_1,T_2}(t_1, t_2)/f_{T_1}(t_1)$. 
Since the numerical algorithm gives the joint density in discretized form it is natural to construct a 
discrete Markov chain with discrete states $x_1, x_2, \ldots , x_n$ and transition matrix
\begin{align}
\Prob(T_2 = x_k \mid T_1 = x_j) 
&= \frac{f_{T_1,T_2}(x_j, x_k)}{\sum_k f_{T_1,T_2}(x_j, x_k)}.
\label{Tmatrix}
\end{align}

Figure~\ref{Fig:ConditionalWH6} illustrates he conditional pdf:s for the shifted Gaussian spectrum 
{\sf WH6} with clearly correlated crossing intervals, $\Corr (T_1,T_2) = 0.40$. 
The figure shows the marginal pdf for a single zero crossing interval and 
four conditional pdf:s for a second interval given the length of the first. 
It is now possible to construct a Markov process with the transition matrix 
\eqref{Tmatrix} and use it as an approximation for the distribution of the 
whole crossing sequence. A natural, exact, and rather strong test of the model 
can be based on the {\it trivariate} distribution of three consecutive intervals. 
The exact tri-variate density  $f_{T_1,T_2,T_3}(t_1, t_2, t_3)$ can be computed by {\sf RIND} 
with the same degree of complexity as the bivariate density. 
 
 \begin{figure}[t!]
\begin{minipage}[c]{0.47\textwidth}
\includegraphics[width=0.8\linewidth]{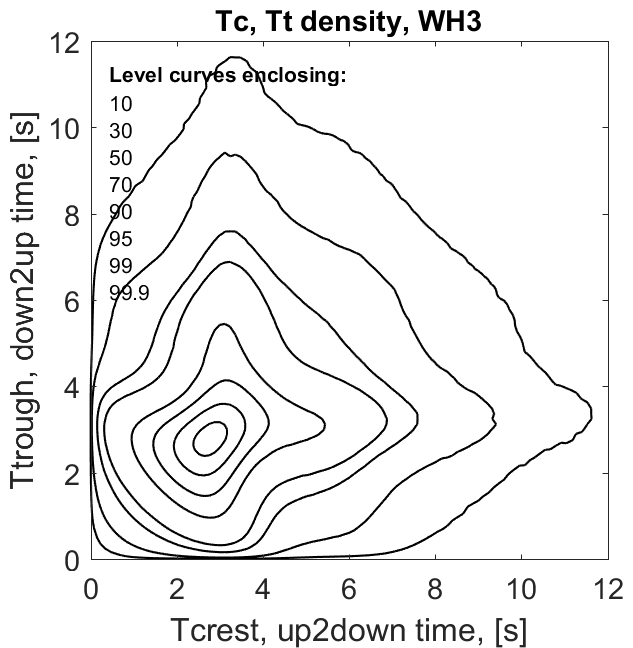}
\end{minipage} 
\begin{minipage}[c]{0.47\textwidth}
\includegraphics[width=0.9\linewidth]{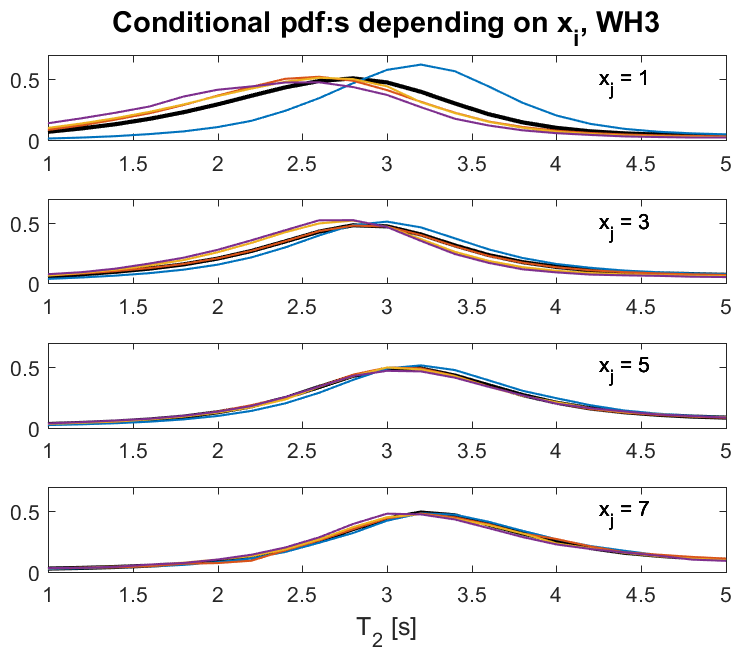}
\end{minipage}
\\
\begin{minipage}[c]{0.47\textwidth}
\includegraphics[width=0.8\linewidth]{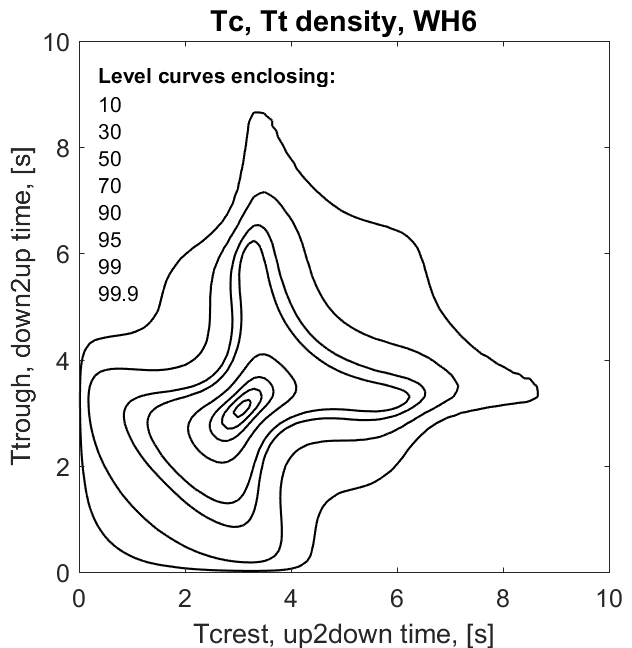}
\end{minipage} 
\begin{minipage}[c]{0.47\textwidth}
\includegraphics[width=0.9\linewidth]{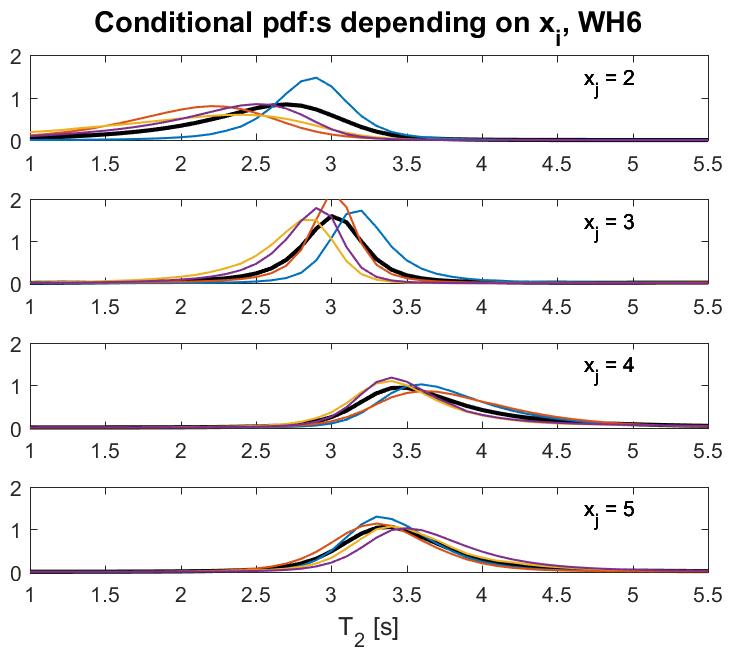}
\end{minipage}
\\
\begin{minipage}[c]{0.47\textwidth}
\includegraphics[width=0.8\linewidth]{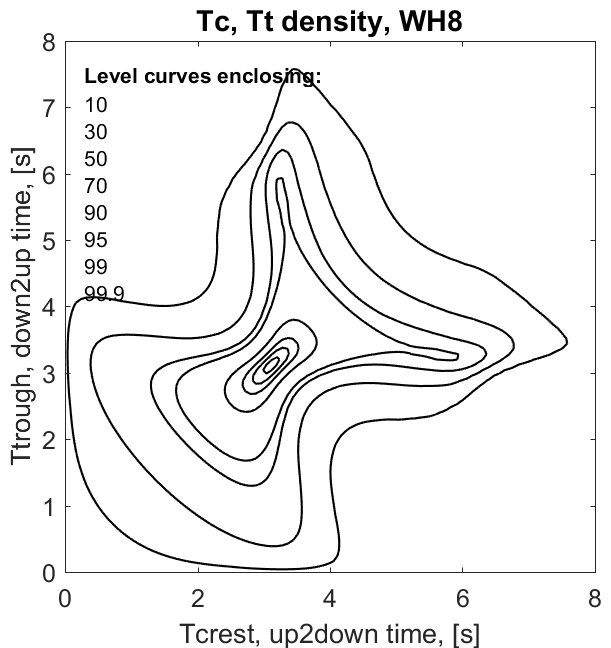}
\end{minipage} 
\begin{minipage}[c]{0.47\textwidth}
\includegraphics[width=0.9\linewidth]{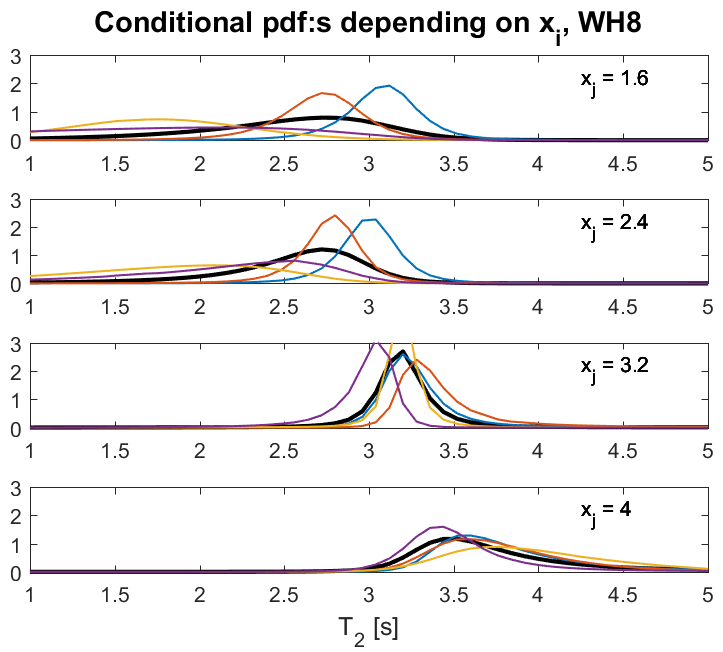}
\end{minipage}
\caption{Markov tests. The top row for {\sf WH3}, the middle for {\sf WH6} and the bottom for {\sl WH8}. 
The left plot shows the 2D density of successive half periods. The right plot shows 
conditional densities for $T_2$, conditioned on the value of $T_1 = x_j$ as black thick lines in each 
subplot. The coloured curves show the effect of also conditioning on the preceding $T_0$ taking any 
of the $x$ values.}
\label{Fig:MarkovtestWH3}
\end{figure}

If the crossing sequence is a Markov chain, then the conditional pdf of an interval $T_2$, given 
the length of the previous interval $T_1$, should be equal to its conditional pdf, given the lengths 
of the two previous ones $(T_0, T_1)$, i.e., for all $x_i, x_j, x_k$, the following equality should hold,
\begin{align}
\Prob (T_2 = x_k \mid T_0 = x_i, T_1 = x_j) &=  
\frac{f_{T_0,T_1,T_2}(x_i, x_j, x_k)}{\sum_{i} f_{T_0,T_1,T_2}(x_j, x_k, x_n)} 
=  \Prob(T_2 = x_k \mid T_1 = x_j). 
\label{Cond3D}
\end{align}

We illustrate the technique on the {\sf WH3, WH6, WH8} processes, as examples of processes 
with ``almost independent'', ``moderately dependent'', and ``strongly dependent'' successive intervals. 

We use the {\sf RIND} applications {\sf cov2ttpdf} and {\sf cov2tttpdf} to compute 
the 2D and 3D densities in \eqref{Tmatrix} and \eqref{Cond3D} and then compute and plot 
the left hand sides conditional densities for selected values of $x_j$ and, for each $x_j$, 
for different $x_i$. Figure~\ref{Fig:MarkovtestWH3} shows the results.

The left plot in each figure shows the 2D density of successive half periods. The right plot shows 
conditional densities for $T_2$, conditioned on the value of $T_1 = x_j$ as black thick lines in each 
subplot. The coloured curves show the effect of also conditioning on the preceding $T_0$ taking any 
of the $x$ values. For example, the blue curve is the conditional pdf, given $T_0 = T_1 = x_1$, etc.  
 
 From the figures, we can conclude that for {\sf WH3} the black curves are different, 
 proving the dependence between intervals. However, the coloured curves in each 
 subplot do not deviate  much from the respective black curve, 
 except for the shortest interval, which indicates that the Markov chain approximation can be used. 
 For {\sf WH6} the deviation from a Markov model is stronger, and for {\sf WH8} the Markov 
 model fails altogether.

\section{Conclusions}
Advances in statistical computing during recent decades has made it possible to compute probabilities and 
expectations in very high-dimensional nearly singular normal distributions. We have illustrated how the 
{\sc Matlab} implementation {\sf RIND} of these methods can be used to solve intricate level crossing problems in 
Gaussian processes. 

We have shown how {\sf RIND} is used to compute the bivariate distribution of the 
distance between three successive level crossings by a general stationary Gaussian process, based only on its 
covariance function. We have identified processes where successive intervals are almost independent, including 
diffusion related processes in different dimensions, where the ``Independent Interval Approximation'' (IIA) is 
expected to provide accurate results due to small correlation between subsequent crossing intervals.  
However even in this favorable situation, the IIA in 2D clearly underestimates ($\approx 0.1863$) the actual value of the persistency exponent (=0.1875), which is defined as the exponential rate with which the tail 
of the interval length distribution falls off. 
An additional problem with the IIA, when it utilizes the covariance matching, is that it does not yield a valid probability distribution.
We conclude that the IIA has to be used with caution even in the cases when the independence is nearly valid. 

On the other hand, we show that our algorithm works remarkably well both for obtaining the actual distribution of the crossing intervals but also in approximating the tail of the length distribution.
We demonstrated it not only for the  nearly independent interval cases, as we have also identified processes with strong dependence  and very complex structure, including standard ocean 
wave models, shifted Gaussian spectra, and simple rational spectra.
For these processes using the IIA is even more questionable, however the  distributions computationally retrieved through {\sf RIND} are matching very well the computationally intensive simulations from the process. 
We have also indicated how the method can 
be extended to deal with higher order dependence and Markov dependence. 

The precision of the algorithm has been illustrated on the distribution of very long crossing intervals. For diffusion in 2D the {\sf RIND} 
algorithm gives the correct theoretical value, and for other dimensions the results agrees with other large Monte Carlo studies. 


\appendix
\section{Appendix}
\subsection{Generalized Rice's formula}

The results in \eqref{A} and \eqref{jointpdf} given in Sections~\ref{sss:distrel} and \ref{ss:DRP} 
can be shown using a generalized Rice's formula presented in 
full generality in \cite[Thm.~6.4]{AzaisW}.  

\begin{theorem}[Weighted Rice Formula] Let $Z : U \rightarrow R^d$ be a random field, $U$ an open
subset of $R^d$, and $u \in R^d$ a fixed point. Assume that:
\begin{itemize}\setlength\itemsep{-1mm}
\item[(i)] $Z$ is Gaussian;
\item[(ii)] almost surely the function $t \rightarrow Z(t)$ is of class $C^2$;
\item[(iii)] for each $t\in U$, the distributtion of $Z(t)$ is non-degenerated, i.e.\ has a 
non-singular covariance matrix;
\end{itemize}
Then for $N_\Lambda(u)$, the number of $t$ in a compact $\Lambda \subset U$ such that $Z(t)=u$,
\begin{eqnarray}
\label{Rice0}
\Ex[N_\Lambda(u)]
=\int_\Lambda \int |\det \dot{z}| f_{\dot{Z}(t),Z(t)}(\dot{z},u)\, \rd\dot{z}\,\rd t.
\end{eqnarray}

In addition, assume that for each $t\in U$ one has another field, $Y^t : 
W\rightarrow  R^n$, defined on some topological space $W$, and verifying the following conditions:

\begin{itemize}\setlength\itemsep{-1mm}
\item[(a)] $Y^t(w)$ is a measurable function of $(\omega, t,w)$ and almost surely, 
$(\omega,t,w)\rightarrow Y^t(w)$ is continuous;
\item[(b)] for each $t\in U$ the random process $(t,w)\rightarrow (Z(t), Y^t(w))$ 
defined on $U\times W$ is Gaussian.
\end{itemize}
Then, if $ g : U\times C(W,R^n) \rightarrow  R$ is a bounded and continuous function,
when $C(W,R^n)$ has the topology of uniform convergence on compact sets, then for
each compact subset $\Lambda$ of U one has
\begin{equation}
\label{Rice1}
\Ex\Bigg[\sum_{t\in \Lambda, Z(t)=u} g(t, Y^t)\Bigg] = \int_\Lambda 
\Ex[|\det \dot{Z}(t)|g(t,Y^t))|Z(t)=u]  f_{Z(t)}(u)\,\rd t.
\end{equation}
\end{theorem}
In the theorem, $Y^t(w)$ is an $n$-dimensional  mark attached to  time point $t$. When using the theorem 
on different crossing problems, we will refer to $n$ as the dimension of the problem. 

There are more general versions of Rice's formula, see \cite{Zahle} and \cite{PodgorskiRM} 
shown under minimal assumptions on a random field. However the price is that they are 
valid for {\it almost} all levels $u$. Consequently the result shown in this appendix are also 
(in practice) applicable for any smooth, symmetrical, time reversible process, 
not necessarily Gaussian. 
We now demonstrate relations \eqref{Aa}, \eqref{Ac}, and \eqref{jointpdf} in this paper. 

Let $X(t)$ be a standardized Gaussian process satisfying the assumptions of the generalized 
Rice's theorem.  
We shall use the theorem to find the distributions of the variables  $A, B$, 
and the ergodic distribution of $T$, the distance between consecutive zeros, and the formula 
for the joint probability density of two consecutive intervals between zeros, $T_1,T_2$. 
Here we let $u=0$, however similar formulas can be shown to hold for any $u$.

Furthermore,  previously used notation will be used, viz. $X_{s,t,u} = (X(s), X(t), X(u))$, 
$\dot{X}_{s,t,u} = (\dot{X}(s), \dot{X}(t), \dot{X}(u))$, 
$\dot{X}_{s,t,u}^{+-+} = \dot{X}(s)^+ \dot{X}(t)^- \dot{X}(u)^+$, 
and $X_{s,t} = (u,v)$ means 
$X(s) = u, X(t) = v$. Moreover, $a \leq X_{s,t} \leq b$ 
means that for each $u \in (s,t): a \leq X(u) \leq b$, 
 while $\{a \leq X_{s,t} \leq b\}$ also stands for the indicator function of this set, 
i.e. the function equal to one whenever the condition between the brackets is 
satisfied and to zero otherwise. 

\subsection{Probability density function of time to first crossing   
and stationary distribution of time between two zeros} 
\subsubsection{Pdf of time $A$ to first zero crossing}
Consider a zero mean stationary smooth Gaussian process $X(t)$ and take $u=0$.  
Let  $Z(t)=X(t)$ and $Y^t(w)=X(w)$ so the dimension is $n=1$.  
(Note that $X$ has a.s. no local extremes with height equal zero.)  
Define  $g$ as the indicator function equal to one if the process has a constant 
sign in the interval $(0,t)$, $g(t,Y^t) = \{X_{0,t} > 0\}+\{X_{0,t} < 0\}$. 

For $\Lambda=[0,a]$ the generalized  Rice's formula writes
\begin{align*}
\Ex\Bigg[\sum_{t\in \Lambda, Z(t)=u} g(t, Y^t)\Bigg] &= \int_0^a \Ex [|\dot{X}(t)|
(\{X_{0,t} > 0\}+\{X_{0,t} < 0\}) \mid X(t)=0] \, f_{X(t)}(0)\,\rd t\\
&=2\int_0^a \Ex [|\dot{X}(t)| \{X_{0,t} > 0\} \mid X(t)=0]\,  f_{X(t)}(0)\,\rd t,
\end{align*} 
since $X$ is symmetrical around $0$. 
Now the weighted sum $\Ex\big[\sum_{t\in \Lambda, Z(t)=0} g(t, Y^t)\big]$ is equal to the 
probability that $A\in\Lambda$, and hence the probability density of $A$ (and $B$) is given by
\begin{align}\label{Durbin0}
f_A(a)&=2\Ex [\dot{X}(a)^-\{X_{0,a} > 0\} \mid X(a)=0]\,  f_{X(a)}(0)\\
&=2\Ex [\dot{X}(0)^-\{X_{-a,0} > 0\} \mid X(0)=0]\,  f_{X(0)}(0)\nonumber\\
&=2\Ex [\dot{X}(0)^+\{X_{0,a} > 0\} \mid X(0)=0]\,  f_{X(0)}(0),\label{Durbin}
\end{align}
by time reversibility of $X$.

Observe that the  function $g(t, Y^t)$ is an indicator function and hence not continuous. 
Hence one needs to show that (\ref{Rice1}) is valid also for such a function $g$. 
This is done by employing the theorem for a sequence of continuous functions $g_\epsilon$, 
which converges to $g$ in such a way that a dominated or monotone convergence theorem can be applied.
For example see the proof of Theorem~7.1 in \cite{AbergRL}. 
A similar approach can be used here to prove the presented results.

\subsubsection{Cdf of inter-crossing distance $T$; proof of \eqref{Ac}}
We turn now to the ergodic (stationary) distribution of $T$, the distance between zeros of $X(t)$. 
Following \cite{LindgrenRychlik1991ISR} we have that the tail probability is equal to 
\[
\Prob^{st}(T>t_0)=\frac{\Ex [\mbox{number of $T^+_i\in [0,1]$ such that $T^+_{i+1}-T^+_{i}>t_0$}]}{\Ex[\mbox{number of $T_i^+\in [0,1]$}]}.
\]
Now for fixed $t_0>0$, let $g(t,Y^t)=\{X_{t,t+t_0} > 0\}+ \{X_{t,t+t_0} < 0\}$ 
\begin{align*}
\Ex \Bigg[\sum_{t\in[0,1] Z(t)=u} g(t, Y^t)\Bigg] &= \int_0^1 \Ex[|\dot{X}(t)|
\{X_{t,t+t_0} > 0\}+ \{X_{t,t+t_0} < 0\} \mid X(t)=0]\,  f_{X(t)}(0)\,\rd t \\
&=2\Ex [\dot{X}(0)^+\{X_{0,t_0} > 0\} \mid X(0)=0]\,  f_{X(0)}(0),
\end{align*} 
by stationarity of $X$. Consequently by (\ref{Durbin})   the ergodic (stationary) 
distribution of $T$ satisfies the following relation
\begin{equation}\label{stationaryT}
\Prob^{st}(T>t_0)=\nu^{-1}\Ex [\dot{X}(0)^+\{X_{0,t_0} > 0\} \mid X(0)=0]\, f_{X(0)}(0)=\mu f_A(t_0),
\end{equation}
where $2\nu = 1/\mu$ is the intensity of zeros. Hence the relation \eqref{Ac} is shown.

The code for evaluation of the distribution of $T$ based on the above result and examples of evaluation for the diffusion in 2D is presented in 
the computational section of this Appendix, see Section~\ref{code}

\subsection{Joint pdf of $A, B$ and pdf of $T$; proof of \eqref{Aa}}\label{app:a3}

In this problem $n=2$, $Z({\bf t})=Z(s,t)=(X(s),X(t))$, $s<0<t$. 
Further $\Lambda=[-b,0] \times [0,a]$ and 
$g({\bf t},Y^{\bf t})= \{X_{s,t} > 0\}+ \{X_{s,t} < 0\}$. Then by Rice's formula
\begin{align*}
&\Ex\Bigg[\sum_{{\bf t}\in \Lambda, Z({\bf t})={\bf 0}} g({\bf t}, Y^{{\bf t}})\Bigg]
= \int_{-b}^0\int_0^a \Ex[|\det \dot{Z}({\bf t})| g({\bf t},Y^{{\bf t}} \mid X_{s,t} = {\bf 0}]  \,
f_{X_{s,t}}({\bf 0})\,\rd t\, \rd s\\
&\hspace{10mm} =2\int_{-b}^0\int_0^a \Ex [\dot{X}(s)^+\dot{X}(t)^-
\{X_{s,t} > 0\} \mid X_{s,t}={\bf 0]}  \,f_{X_{s,t}}({\bf 0})\, \rd t\, \rd s.
\end{align*} 
Obviously  the weighted sum $\Ex \big[\sum_{{\bf t}\in \Lambda, Z({\bf t})=0} g({\bf t}, Y^{{\bf t}})\big]$ 
is equal to  $\Prob( B\le b, A\le a)$.  Hence the joint probability density of $B, A$ is given by
\begin{align}
f_{B, A}(b,a)&=2 \Ex[\dot{X}(-b)^+\dot{X}(a)^-
\{X_{-b,a} > 0\} \mid X_{-b,a}={\bf 0}] \, f_{X_{-b,a}}({\bf 0}) \nonumber \\
&=2\,\Ex[\dot{X}(0)^+\dot{X}(a+b)^-
\{X_{0,a+b} > 0\} \mid X_{0,a+b}={\bf 0}] \, f_{X_{0,a+b}}({\bf 0}),  \label{jointBA}
\end{align}
 $0<a,b$, by stationarity of $X$.

Next we  give the pdf of the distribution in (\ref{stationaryT}). 
Again $n=2$, but with $Z({\bf t})=Z(t,s)=(X(t),X(t+s))$, 
$\Lambda=[0,1] \times [0,t_0]$ and $g(t,Y^{\bf t})=\{X_{t,t+s} > 0\}+ \{X_{t,t+s} < 0\}$. 
Obviously $\det \dot{Z}({\bf t})=\dot{X}(t)\dot{X}(t+s)$ and by Rice's formula
\begin{align*}
&E\Bigg[\sum_{{\bf t}\in \Lambda, Z({\bf t})={\bf 0}} g({\bf t}, Y^{{\bf t}})\Bigg]  
= \int_0^1\int_0^{t_0} \Ex[|\det \dot{Z}({\bf t})|g({\bf t}, Y^{{\bf t}}) \mid X_{t,t+s}={\bf 0}]  \,
f_{X_{t,t+s}}({\bf 0})\, \rd t\,\rd s \\
&\hspace{10mm} = 2\int_0^{t_0} \Ex[\dot{X}(0)^+\dot{X}(s)^-
\{X_{0,s} > 0\} \mid X_{0,s}={\bf 0}]  \,f_{X_{0,s}}({\bf 0})\, \rd s,
\end{align*} 
by stationarity.
Consequently the pdf of $T$ is given by
\begin{equation}\label{stationary2}
f_{T}(t)=\nu^{-1} \Ex[\dot{X}(0)^+\dot{X}(t)^-
\{X_{0,t} > 0\} \mid X_{0,t}={\bf 0}] \, f_{{X_0,t}}({\bf 0}).
\end{equation}
A comparison with  \eqref{jointBA} gives 
\[
f_{B, A}(b,a)=\mu^{-1}\,f_T(a+b), \quad 0<a,b,
\]
and  \eqref{Aa} is proved, since $2\nu = 1/\mu$. 
Note that in both the presented examples the probability density function 
of $Z(t,s)$ is not bounded for all $(s,t)\in \Lambda$. 
This is only a technical  problem which is solved by replacing the interval  $[0,t_0]$ in the
definitions of $\Lambda$ by $[\varepsilon, t_0]$, $\varepsilon>0$ and then letting $\varepsilon$ tend to zero.  

\subsection{Stationary joint probability density function of $(T_1,T_2)$}\label{jointT1T2proof}

Now $n=3$ and ${\bf t}=(r,t,s)$, $Z({\bf t})=(X(t+r),X(t),X(t+s))$.
 Furthermore we let $Y^{\bf t}(\tau)=X(\tau)$
and  $\Lambda=[t_1,0]\times [0,1] \times [0,t_2]$ and for $r<0<s$
\[
g({\bf t},Y^{\bf t})=\{X_{t+r,t} <0\}\{X_{t,t+s} > 0\}+\{X_{t+r,t} >0\}\{X_{t,t+s} < 0\}.
\]
Since $\det \dot{Z}({\bf t})=\dot{X}(t+r)\dot{X}(t)\dot{X}(t+s)$ and by symmetry of $X$ around level zero,
 Rice's formula gives 
\begin{align*}
\Ex\Bigg[&\sum_{{\bf t}\in \Lambda, Z({\bf t})=0} \,g({\bf t}, Y^{{\bf t}})\Bigg] 
=\int_{t_1}^0 \int_0^1\int_0^{t_2} \Ex[|\det \dot{Z}({\bf t})|
g({\bf t},Y^{\bf t}) \mid Z({\bf t})={\bf 0}]\,  f_{Z({\bf t})}({\bf 0})\,\rd r\,\rd t\,\rd s\\
&= 2\int_{t_1}^0\int_0^{t_2} \Ex [\dot{X}(r)^-\dot{X}(0)^+\dot{X}(s)^-\{X_{r,0} <0\}
\{X_{0,s} > 0\} \mid X_{r,0,s}={\bf 0}]\, f_{X_{r,0,s}}({\bf 0})\,\rd r \, \rd s,
\end{align*} 
by stationarity and hence the joint probability density function of $T_1<0<T_2$  is  
\begin{equation*}\label{stationary3}
f_{T_1,T_2}(t_1,t_2)=\nu^{-1}\,\Ex [\dot{X}(t_1)^-\dot{X}(0)^+\dot{X}(t_2)^-\{X_{t_1,0} <0\}
\{X_{0,t_2} > 0\} \mid X_{t_1,0,t_2}=({\bf 0})] \,
f_{X_{t_1,0,t_2}}({\bf 0}),
\end{equation*}
in agreement with \eqref{jointpdf}.

\subsection{The validity of the IIA}
\label{Val-IIA}
In this appendix we consider the inverse problem, that constitutes a basis for the IIA approach. 
Namely,  for a given covariance function $R$, normalized so that $R(0)=1$, we ask whether the 
function $\Psi$ defined through  
\begin{equation}
\label{InvEq}
\Psi(s)=\frac{2-s\mu(1-s\mathcal L R(s))}{2+s\mu(1-s\mathcal L R(s))}
=\frac{2}{1+\frac{\mu }2(s-s^2\mathcal L R(s))}-1,
\end{equation}
for a certain  $\mu>0$ corresponds to the Laplace transform of a probability 
distribution function of a non-negative random variable. 

We recall that by Bernstein's theorem the function $\Psi$ is a Laplace transform of a 
probability distribution on the positive half-line if and only if $\Psi(0)=1$, it has all derivatives on $(0,\infty)$ and
$$
(-1)^n \Psi^{(n)}(s)\ge 0,~~s>0. 
$$
It is not difficult to verify that the function given in (\ref{InvEq}) takes value one at $s=0$ 
and it satisfies the above for $n=0$. 
However, as we will see next, covariances $R$ that lead to a valid probability distributions are 
restricted only to ones that have a singular component in their spectral measure, 
the requirement that is rarely satisfied for covariances of interest.
In particular for  Gaussian processes considered in this work,  $\Psi$ is not the Laplace transform of a probability distribution.  
We first formulate this result and then provide the argument through some auxiliary 
facts of a more general nature.  
\begin{proposition}
\label{Gauss}
Let $X$ be a Gaussian process with the (symmetric) spectrum $S$, i.e.\ its covariance is given by
$$
r(t)=\int_{-\infty}^\infty e^{i\omega t} S(\omega)\, \rd \omega.
$$
Then for the covariance $R=R_0^c$ of the clipped process $D_c(t)=\sgn(X(t))$ that has the form
$$
R(t)=\frac{2}{\pi} \arcsin (r(t)/r(0)),
$$
the function  $\Psi$ given through (\ref{InvEq}) does not represent a valid probability 
distribution for any choice of $\mu >0$. 
\end{proposition}
\begin{remark}
We note that  most of physically interpretable processes are given by a symmetric spectrum, 
including these in Table~\ref{Tab:1}. 
In particular, the diffusion process in dimension two has the explicit spectrum 
$S(\omega)= \operatorname{sech}(\pi \omega)$. 
\end{remark}

The argument that shows the above conclusions is based on discussing when the 
Bernstein condition for $n=1$ can be satisfied by $\Psi$, i.e.\ for $s>0$ we investigate if 
\begin{align*}
-\Psi'(s)&=\mu\frac{1-2s\mathcal L R(s)-s^2\mathcal L R'(s)}%
{\left(1+\frac{\mu }2s^2(1/s-\mathcal L R(s))\right)^2}\ge 0.
\end{align*}
By Bochner's theorem any covariance $R$ can be written in  terms of its spectral measure $S_R$ as 
$$
\int_{-\infty}^\infty e^{it\omega}\, \rd S_R(\omega).
$$
From this we have
\begin{align*}
\mathcal L R(s)&=2\int_{0}^\infty  \frac{s}{s^2+\omega^2}\, \rd S_R(\omega),\\
\mathcal L R'(s)&=2\int_{0}^\infty  \frac{\omega^2 - s^2}{s^2+\omega^2}\, \rd S_R(\omega),
\end{align*}
and, by straightforward calculation,
the necessary condition for a covariance $R(t)=\int_{-\infty}^\infty e^{i\omega t}\, \rd S_R(\omega)$ 
to be obtained as a covariance of a switching process is
\begin{equation}
\labmarg{spin}
\int_{0}^\infty \frac{(\omega/s)^4}{\left(1+(\omega/s)^2\right)^2}\, \rd S_R(\omega) \ge
\int_{0}^\infty \frac{1+4(\omega/s)^2}{\left(1+(\omega/s)^2\right)^2}\, \rd S_R(\omega),
\end{equation}
an inequality that has to be satisfied for each $s>0$.

Let us consider now that $\rd S_R(\omega)=S_R(\omega)\, \rd \omega$. 
Since the normalized spectrum $S_R$ integrates to one (variance of the symmetric switching 
process is always one), there exists $K>0$ such that for all $s>K$ we have
$$
\int_0^{s} S_R(\omega)\, \rd \omega > 4\int_{s}^\infty S_R(\omega)\, \rd \omega.
$$
Then for $s>K$ we have
\begin{align*}
\int_{0}^\infty \frac{u^4-1-4u^2}{(1+u^2)^2} S_R(us)\, \rd u &=
\int_{0}^{1} \frac{u^4-1-4u^2}{(1+u^2)^2} S_R(us)\, \rd u+\int_{1}^\infty 
\frac{u^4-1-4u^2}{(1+u^2)^2} S_R(us)\, \rd u\\
&
 \le - \int_{0}^{1} \frac{ S_R(us)}{(1+u^2)^2}\, \rd u + \int_{1}^\infty  S_R(us)\, \rd u,
\end{align*} 
where the last integral appears because $(u^4-1-4u^2)/(1+u^2)^2<1$.  
Further, 
\begin{align*}
\int_{0}^\infty \frac{u^4-1-4u^2}{(1+u^2)^2} S_R(us)\, \rd u
 &
 \le -\frac 1 4\int_{0}^{1}  S_R(us)\, \rd u + \int_{1}^\infty  S_R(us)\, \rd u\\
 &= \frac{1}{4s}\left(4\int_{s}^\infty  S_R(\omega)\, \rd \omega 
 - \int_{0}^{s}  S_R(\omega)\, \rd \omega\right) <0. 
\end{align*}

This proves the following result.

\begin{proposition}  If the covariance $R$ in (\ref{InvEq}) is given by a spectrum through 
$$
\int_{-\infty}^\infty e^{it\omega} S_R(\omega)\, \rd \omega,
$$
 then the formula does not yield $\Psi$ corresponding to a probability distribution for any $\mu>0$.
 \end{proposition} 
To obtain Proposition~\ref{Gauss} from the above arguments 
it is enough to show that the covariance given through 
$$
R(t)=\frac{2}{\pi} \arcsin (r(t)/r(0)),
$$
has a continuous spectrum $S_R$ whenever $r$ has such. 
This follows from the series expansion 
$$
\operatorname{arcsin}z=z+\frac{1}{2}\frac{z^{3}}{3}+\frac{1\cdot 3}{2\cdot 4}%
\frac{z^{5}}{5}+\frac{1\cdot 3\cdot 5}{2\cdot 4\cdot 6}\frac{z^{7}}{7}+\cdots,
$$
that implies that for the (normalized) spectrum $S_r$ of $r(t)/r(0)$ the spectrum of $R$ is given by
$$
S_R=\frac{2}{\pi}\left( S_r+\frac{1}{2}\frac{S_r^{*3}}{3}+\frac{1\cdot 3}{2\cdot 4}\frac{S_r^{*5}}{5}
+\frac{1\cdot 3\cdot 5}{2\cdot 4\cdot 6}\frac{S_r^{*^7}}{7}+\cdots\right).
$$

\subsubsection{Diffusion in dimension two}
\label{ddt}
We illustrate the problem of recovering the crossing interval distribution by the IIA approach through 
the analysis of the diffusion in dimension two, with covariance as given in  {\sf BMS2}. 
This example can serve as the benchmark case since a lot of computation can be made 
more explicit and the actual value  of persistency  $\theta(2)=3/16=0.1875$ has been recently established in \cite{PoplavskyiSchehr2018}.

First, we report the following Laplace transform of the covariance 
obtained through a series expansion of $\displaystyle \arcsin x$,
\begin{multline*}
\mathcal L\left(\frac 2 \pi \arcsin\left(\frac{2}{e^{t/2}+e^{-t/2}}\right)\right)(s)
=
\mbox{\large $\displaystyle \frac 2 \pi \sum_{l=0}^{\infty}$}~
\frac{\displaystyle \sum_{k=0}^{l}
\frac{1}{
2^{k}(2k+1)}
\binom{l+k}{k}\binom{s-1/2+k}{k}
}{\displaystyle 2^l(l+1) \binom{s+1/2+l}{l+1}},
\end{multline*}
where 
$$
\binom{x}{y}=\frac{\Gamma(x+1)}{\Gamma(x-y+1)\Gamma(y+1)}.
$$
We note the following structure: 
\begin{align*}
\mathcal L\left(\frac 2 \pi \arcsin\left(\frac{2}{e^{t/2}+e^{-t/2}}\right)\right)(s) 
=\sum_{l=0}^\infty \frac{P_l(s)}{(s+1/2)\dots (s+1/2+l)},
\end{align*}
where $P_l$ are polynomials in $s$ of the order $l$,  given by
\begin{align*}
P_l(s)
=\frac{l! }{\pi 2^{l-1}}\left(
1+
\sum_{k=1}^{l}
\frac{\binom{l+k}{k}
}{
2^{k}(2k+1)k!}  
{(s+1/2)\dots (s+1/2 +k-1)}
\right).
\end{align*}
Consequently, any partial sum that approximates the Laplace transform,
$$
\mathcal L_L(s)=\sum_{l=0}^{L} \frac{P_l(s)}{(s+1/2)\dots (s+1/2+l)}
= \frac{Q_L(s)}{(s+1/2)\dots (s+1/2+L)},
$$ 
 is also a rational function with the numerator being a polynomial in $s$ of order $L$, say $Q_L(s)$, 
 while the denominator is a factorized polynomial of order $L+1$. 
This leads to the following approximate formula for the Laplace transform of the 
distribution of the crossing intervals,
\begin{align}
\Psi_L(s)
\label{Lapp}
&=\frac{4
(s+1/2)\dots (s+1/2+L)
}{
(2+s\mu)(s+1/2)\dots (s+1/2+L)-
  \mu s^2 Q_L(s)
}-1.
\end{align}
We note that the first term is a ratio of two polynomials, with the numerator of order $L$ and 
denominator of order $L+1$.
We immediately conclude that  this cannot be a Laplace transform of a probability distribution 
since the function equal to minus one is the Laplace transform of a negative atomic measure 
with atom minus one at zero,  while the ratio of the polynomials with the order of the numerator smaller than that of the denominator is incapable of canceling 
this negative atomic measure. 

For the diffusion in D2, the recommended value of $\mu$ is $2\pi$ which is the 
average value of the zero crossing interval. 
We mark here that since the actual measure corresponding 
to $\Psi(s)$ is not probabilistic, this choice of $\mu$ is only partially justified. 

\begin{figure}[t!]
\includegraphics[width=\textwidth]{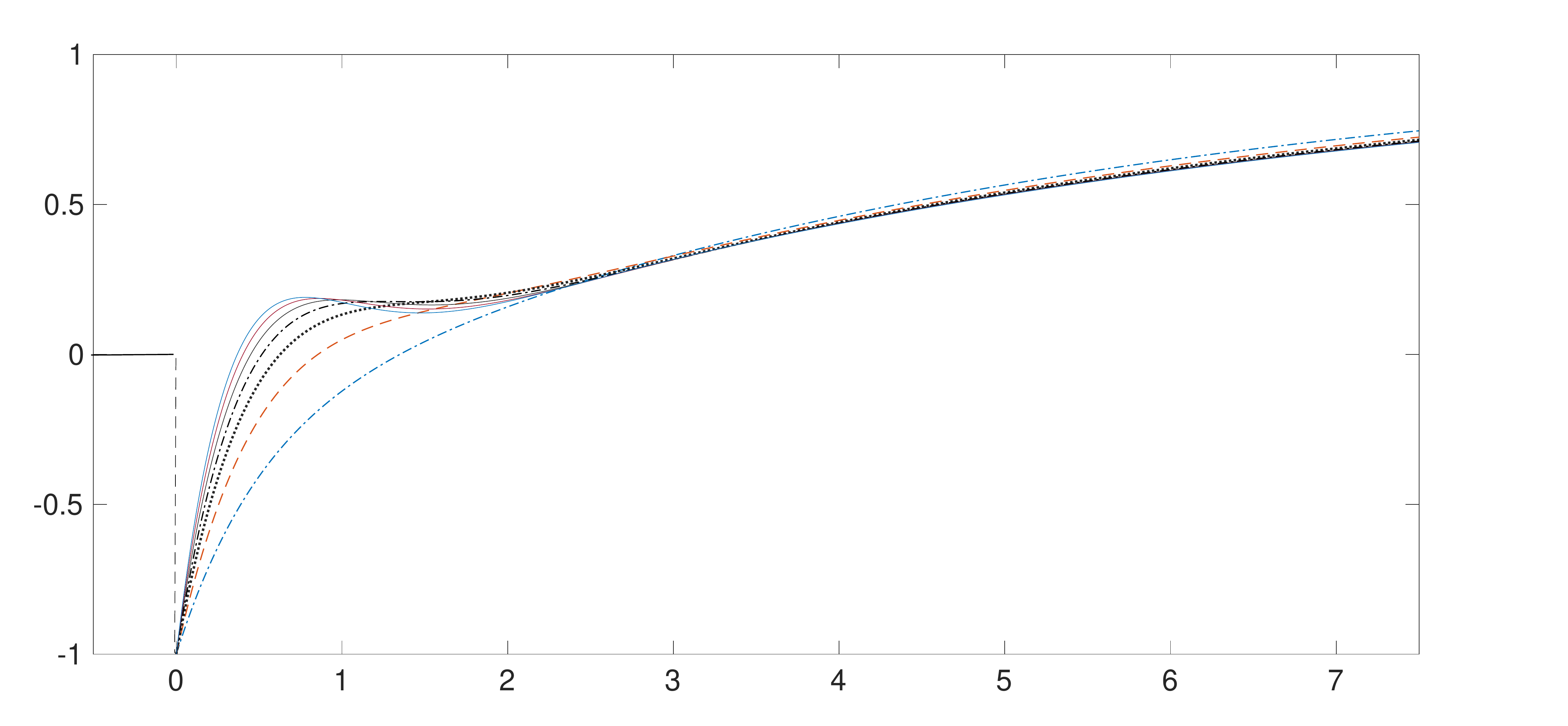}
\caption{Several cumulative non-probabilistic distributions representing signed measures 
obtained through approximation (\ref{Lapp}) for $L=0,1,2,3,4,5,6$.
The higher order approximations attempt to correct the high `bump' near the origin.}
\label{app}
\end{figure}

Consider, for example, $L=0$. 
Direct computations leads to the following zero order approximation
\begin{align*}
\Psi_0(s)
&=
\left(
\frac{
a
}{s-s_1}
+\frac{
b
}{s-s_2}
\right)-1,
\end{align*}
where 
\begin{align*} 
a&\approx 0.2740,~~b\approx 1.4779\\
s_1&\approx -0.2150,~~s_2 \approx -2.0369
\end{align*}
and the corresponding ``distribution'' that is obtained from the inverse Laplace transform 
has the absolute continuous part given by the density
\begin{align*}
f(t)=a e^{t s_1} + be^{ts_2}.
\end{align*}
We point out that  the inverse Laplace method does not produce a probabilistic measure due 
to the atom at zero although the total mass is still equal to one.
The equivalent to the expectation is equal to $\mu=2\pi$ but it cannot be interpreted as the 
average value of the crossing interval. 
The ``cdf'' of this measure is presented in Figure~\ref{app}. 
Finally, we note that the asymptotics is governed by the  exponent $-s_1\approx 0.2150$, 
which could be viewed as a crude ($L=0$) approximation of the persistency exponent. 

Approximations for other values of $L$ can be done similarly. 
The obtained graphs of the quasi cumulative distribution functions are shown in Figure~\ref{app}. 
It can be seen that there is a probabilistically uninterpretable negative jump at zero, 
which higher order approximations attempt to `correct'.
This however leads to a probabilistically uninterpretable loss of the monotonicity on the 
continuous part for small values of the time crossings.
However, as previously reported in the literature,  the approach gives a reasonable estimate of the 
tail of the distribution as can be seen in the obtained persistency coefficients. 
In Table~\ref{perco}, we show the approximated values of the persistency coefficient.
Recall that by a recent result  in \cite{PoplavskyiSchehr2018}, the true persistency for the diffusion in two 
dimensions is $\alpha=3/16=0.1875$. 
We observe that this value very closely followed by even low order approximations and 
is actually attained by the approximation of order six. 
However, further approximations are producing smaller values eventually stabilizing at 
$\theta_{IIA}=0.1863$, which is consistent with the results reported in 
\cite{BrayMajumdarSchehrAiP2013} and \cite{MajumdarSBC}. 
Our conjecture is that this is the actual value of the persistency obtained by the IIA 
approach and thus it underestimates the actual value of the parameter. 

\begin{table}
\caption{The approximated values of the persistency $\theta_L$ using the largest 
negative poles of $\Psi_L$ given in (\ref{Lapp}).  }
\label{perco}
\begin{center}
\begin{tabular}{ccccccccc}
\toprule
$L$ & 0 & 1 & 2 & 3 & 4 & 5 & 6 & 7   \\
\midrule 
$\theta_L$ & 0.2150 & 0.1991 & 0.1930 & 0.1902 & 0.1887 & 0.1880 & 0.1875 & 0.1872  \\
\bottomrule
\toprule
$L$ & 8 & 9& 10 & 15 & 20 & 30 & 50 &80  \\
\midrule 
$\theta_L$ & 0.1870  &  0.1869 & 0.1868 & 0.1866 & 0.1865 & 0.1863 &  0.1863  &   0.1863\\
\bottomrule
\end{tabular}
\end{center}
\end{table}

\subsection{Examples of the code for some numerical results}
\label{code}
In this section, for the reader and potential user convenience, we provide a complete 
code of some numerical evaluation presented in the paper using the {\sf WAFO}-package.  

\begin{figure}[t!]
\centerline{
\includegraphics[width=0.95\textwidth]{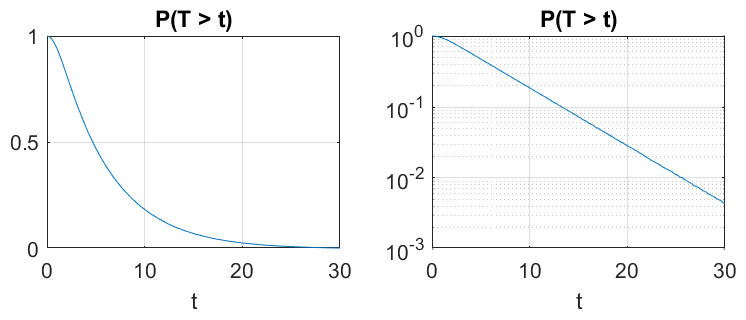}
}
\caption{$P(T>t)$ for diffusion {\sf BMS2} computed by {\sf RIND}, in the regular scale {\it (Left)} and in the logarithmic scale {\it (Right)}.}
\label{fig:PTgtt}
\end{figure}

The first code shows computation of the distribution of the time $T$  between subsequent crossings. 
In Figure~\ref{fig:PTgtt}, the resulting cdf in the normal and logarithmic scales is presented. 

\begin{lstlisting}[language=Matlab]
% Script to compute the persistence, 1-cdf of intercrossing time
dt=0.2; N=; t=0:dt:N*dt;
Rt=sech(t/2);
Rt1=tanh(t/2).*sech(t/2)/2;
Rt2=(tanh(t(1)).*Rt1(1)+Rt(1).*(1-tanh(t(1)/2).^2)/2)/2;

R=[toeplitz(Rt(1:end-1)); fliplr(Rt1(2:end))]; 
R=[R [R(end,1:end) Rt2]'];
R=[R;[fliplr(Rt(2:end)) 0]]; 
R=[R [R(end,1:end) Rt(1)]'];

% Set lower and upper bounds and compute P, prob to stay between bounds
u=0; Blo=[u 0]; Bup=[inf inf];
Pt=[];
for i=1:N
    Nt=N-i+1; xc=u; indI=[0 Nt Nt+1]; mi = zeros(Nt+2,1);
    P=rind(R(i:end,i:end),mi,Blo,Bup,indI,xc,Nt);
    Pt= [Pt; [t(Nt+1) P]];
end
mu=sqrt(Rt2/Rt(1))/2/pi;
Pt(:,2)=Pt(:,2)/mu;
Pt=[Pt;[0 1]];
\end{lstlisting}\vspace{3mm}

In the next block of Matlab-code,  some other computational results used across the paper are presented to illustrate relative simplicity of utilizing {\sf RIND}. For description of the special command {\sf rindopset} see the help text in {\sf WAFO}. 

{\footnotesize
\begin{lstlisting}[language=Matlab]
function [f,find] = cov2ttpdf_cov(cc,level,paramt,modify,varargin)
%COV2TTPDF computes bivariate (Tc,Tt) 
% 
%  CALL [f,find] = cov2ttpdf_cov(c,level,paramt,modify,options)
%
%   f       = 2D pdf structure for pair of crest/trough period
%   find    = 2D pdf structure with assumed independent marginals
%
%   c       = symbolic covariance function
%   level   = reference level, default = 0
%   paramt  = [Tc Td Nc Nd] = [real real integer integer] 
%               max crest/trough period and subdivision 
%               with Nc/Nt = Tc/Td (if not, Tc,Td will be adjusted) 
%               Default: Tc = Td = 10, Nc = Nd = 20 gives
%				pdf at x1,x2 = 0:0.5:10
%   modify  = 1, modify irregular case (when L3 neq 0)
%             0, do not modify (default)
%   options = options for RIND, set by  rindoptset
%
% Use spec2ttsim to simulate half and full periods  

% History: 
% Updated 2018 and 2019 to allow irregular spectra
% Revised by GL May 2017 - the input c (covariance function 
%   or spectrum structure) is now normalized automatically to 
%   m0 = m2 = 1
% Made by Georg Lindgren, February 2017, 
% based on IR, PAB, KP and others spec2tpdf WAFO routine
% Used on Matlab 2017b

start = clock;
narginchk(1,inf)

defaultSpeed = 2;
defaultMethod = 5;
defaultoptions = rindoptset('speed',defaultSpeed,'method',defaultMethod);
if nargin<5
	opt = defaultoptions;
else
    opt = rindoptset(defaultoptions,varargin{:});
end
    
if nargin < 4 || isempty(modify)
    modify = false;
elseif modify == 1
    modify = true;
else
    modify = false;
end
if nargin<3 || isempty(paramt)
    paramt = [10 10 20 20];
end
if  nargin < 2 || isempty(level)
    u = 0;
else 
    u = level;
end

T1 = paramt(1); % Aimed max crest interval
T2 = paramt(2); % Aimed max trough interval
Nc = paramt(3); % Nc+1 = length(x1)
Nt = paramt(4); % Nt+1 = length(x2)
Ttot = T1+T2;
Ntot = Nc+Nt;
TimeStep = Ttot/Ntot;
T1 = Nc*TimeStep; % Final max crest interval
T2 = Nt*TimeStep; % Final max trough interval

c = cc;
covtxt = formula(c);

%%%%%%%%%%%%%%%%%%%%%%%%%%%%%%%%%%%%%%%%%%%%%%%%%%%%%%
% Create pdf-structure and prepare for information
    f = createpdf(2);
    f.labx{1} = 'Tcrest, up2down time, [s]';
    f.labx{2} = 'Ttrough, down2up time, [s]';
    f.title = 'Tc, Tt density';
    f.level = level;
    f.note = ['cov = ' char(covtxt)];
    f.date = start;
    f.x{1} = linspace(0,T1,Nc+1);
    f.x{2} = linspace(0,T2,Nt+1);
    f.f = zeros(Nc+1,Nt+1);
    time = linspace(0,T1+T2,Nc+Nt+1);
% End: Create pdf-structure
    
%%%%%%%%%%%%%%%%%%%%%%%%%%%%%%%%%%%%%%%%%%%%%%%%%%%%%%
% Compute R = covariance function and two derivatives 
% from symbolic covariance function and normalize to L0 = L2 = 1
% Note that irregular case is allowed
    syms t normedt
    L2_ = diff(c,2);
    L2_ = limit(L2_,0);
    korr = sqrt(-L2_);
    normedt = t/korr;
    c(t) = cc(normedt);

    ttp=time(time>0);
    dt = ttp(2)-ttp(1);
    t=sym(ttp);
    R.t=time;
    y=c(t);
    y0=limit(c,0);
    R.R=[double(y0) double(y)];
    L0=double(y0);

    dc=diff(c,1);
    y=dc(t);
    y0=limit(dc,0);
    R.Rt=[double(y0) double(y)];
    L1=double(y0);

    ddc=diff(c,2);
    y=ddc(t);
    y0=limit(ddc,0);
    R.Rtt=[double(y0) double(y)];
    L2 = -double(y0);
    
    dddc=diff(c,3);
    syms t
    L3=double(limit(dddc,t,0,'right')); 
    %Case is irregular if  L3 not equal 0
    
    ddddc=diff(c,4);
    L4=double(limit(ddddc,t,0,'right'));
% End: Covariances from symbolic covariance function

%%%%%%%%%%%%%%%%%%%%%%%%%%%%%%%%%%%%%%%%%%%%%%%%%%%%%%%
% Compute all covariance matrices for process and derivative
    R11 = R.R(1:end);
    R11 = toeplitz(R11);
    R22 = -R.Rtt(1:end);
    R22 = toeplitz(R22);
    R12 = toeplitz(R.Rt(1:end));
    R12 = diag(diag(R12)) + triu(R12,1) - tril(R12,-1);
    R21 = R12';
% End: Compute all covariance matrices

%%%%%%%%%%%%%%%%%%%%%%%%%%%%%%%%%%%%%%%%%%%%%%%%%%%%%%%
% Prepare for RIND
    Xc = zeros(3,1); 
    Mine=100; ERR = zeros(Nc,Nt); TERR = zeros(Nc,Nt);
    tstart = tic;
% End: Prepare for RIND

%%%%%%%%%%%%%%%%%%%%%%%%%%%%%%%%%%%%%%%%%%%%%%%%%%%%%%%
% Run RIND for all point in 2D grid
for j=1:Nc
    waitbar(j/Nc)
    for k = 1:Nt
    % Select sub-covariance matrices from RRR
    Btt = [R11(j+1:Nc,j+1:Nc)         -R11(j+1:Nc,Nc+2:Nc+Nt+1-k);
          -R11(Nc+2:Nc+Nt+1-k,j+1:Nc)  R11(Nc+2:Nc+Nt+1-k,Nc+2:Nc+Nt+1-k)];
          
    Btd = [R12(j+1:Nc,j) -R12(j+1:Nc,Nc+1) R12(j+1:Nc,Nc+Nt+2-k) ;
          -R12(Nc+2:Nc+Nt+1-k,j)  R12(Nc+2:Nc+Nt+1-k,Nc+1) ...
                 -R12(Nc+2:Nc+Nt+1-k,Nc+Nt+2-k)];
          
    Btc = [R11(j+1:Nc,j)  R11(j+1:Nc,Nc+1)  R11(j+1:Nc,Nc+Nt+2-k);
          -R11(Nc+2:Nc+Nt+1-k,j) -R11(Nc+2:Nc+Nt+1-k,Nc+1) ...
                  -R11(Nc+2:Nc+Nt+1-k,Nc+Nt+2-k)];
          
    Bdd = [R22(j,j)         -R22(j,Nc+1)          R22(j,Nc+Nt+2-k);
          -R22(Nc+1,j)       R22(Nc+1,Nc+1)      -R22(Nc+1,Nc+Nt+2-k);
          R22(Nc+Nt+2-k,j) -R22(Nc+Nt+2-k,Nc+1)  R22(Nc+Nt+2-k,Nc+Nt+2-k)];
           
    Bdc = [R21(j,j)         R21(j,Nc+1)         R21(j,Nc+Nt+2-k);
          -R21(Nc+1,j)     -R21(Nc+1,Nc+1)     -R21(Nc+1,Nc+Nt+2-k);
           R21(Nc+Nt+2-k,j) R21(Nc+Nt+2-k,Nc+1) R21(Nc+Nt+2-k,Nc+Nt+2-k)];
           
    Bcc = [R11(j,j)          R11(j,Nc+1)         R11(j,Nc+Nt+2-k);
           R11(Nc+1,j)       R11(Nc+1,Nc+1)      R11(Nc+1,Nc+Nt+2-k);
           R11(Nc+Nt+2-k,j)  R11(Nc+Nt+2-k,Nc+1) R11(Nc+Nt+2-k,Nc+Nt+2-k)];
           
        BBB = [Btt  Btd  Btc; 
               Btd' Bdd  Bdc; 
               Btc' Bdc' Bcc];
        BBB = BBB + 0.0000001*eye(size(BBB)); % To make sure
        Mine = min(Mine,min(eig(BBB))); % A possible check
   % End: Select sub-covariance       
        
   % Define gridpoint specific input for RIND
        Ntdc = length(BBB);
        Ntimes = length(Btt); Nt1 = int8(Nc-j); Nt2 = int8(Nt-k);
        Nder = length(Bdd);
        m = [repmat(-u,Nt1,1); repmat(u,Nt2,1); [0 0 0 -u -u -u]'];
        Indb = [0:Ntimes+Nder];
        Blo = zeros(1,Ntimes+Nder);
        Bup = Inf(1,Ntimes+Nder);
   % End: Define grigpoint specific input 
        
   % Call RIND
        [ff,err,terr] = rind(BBB,m,Blo,Bup,Indb,Xc,Ntimes,opt);
        f.f(Nc+2-j,Nt+2-k) = ff*2*pi*exp(u^2/2);
        ERR(j,k) = err;
        TERR(j,k) = terr;
   end
end
% RIND completed

%%%%%%%%%%%%%%%%%%%%%%%%%%%%%%%%%%%%%%%%%%%%%%%%%%%%%%
% Sum up for 2D output pdf  f
telapsed = toc(tstart);
f.time = telapsed;
f.f = f.f';
if u==0
	f.f = (f.f+f.f')/2; % If you are sure it should be symmetric
end

f.integral = trapz(f.x{1},trapz(f.x{2},f.f,1));
f.opt = opt;
f.marginal1 = trapz(f.x{2},f.f);
f.marginal2 = trapz(f.x{1},f.f');
[ql,PL] = qlevels(f.f(:),[10 30 50 70 90 95 99 99.9000]);
f.cl = ql;
f.pl = PL;
f.ERR = ERR;
f.TERR = TERR;

fm1 = trapz(f.x{1},f.x{1}.*f.marginal1);
fm2 = trapz(f.x{2},f.x{2}.*f.marginal2);
fm11 = trapz(f.x{1},f.x{1}.^2.*f.marginal1);
fm22 = trapz(f.x{2},f.x{2}.^2.*f.marginal2);
fm12 = trapz(f.x{1},trapz(f.x{2},(f.x{2}'*f.x{1}).*f.f,1));
f.corr = (fm12-fm1*fm2)/sqrt((fm11-fm1^2)*(fm22-fm2^2));
f.mean = [fm1 fm2];
% End: Sum up output pdf  f

%%%%%%%%%%%%%%%%%%%%%%%%%%%%%%%%%%%%%%%%%%%%%%%%%%%%%%
% Compute pdf  find as product of marginal pdf's and find deviation
if nargout >1
    find = f;
    find = rmfield(find,'corr');
    find.f = f.marginal2'*f.marginal1;
    find.integral = trapz(find.x{1},trapz(find.x{2},find.f,1));
    find.note2 = 'Assumed independence';

    rati = (f.f+eps)./(find.f+eps)/sum(sum(f.f+eps))*sum(sum(find.f+eps));
    num = ~isnan(rati);
    KL = sum(sum((f.f(num)+eps).*log(rati(num))));
    f.KL = KL*((f.x{1}(2)-f.x{1}(1))*(f.x{2}(2)-f.x{2}(1)));
    [ql,PL] = qlevels(find.f(:),[10 30 50 70 90 95 99 99.9000]);
    find.cl = ql;
    find.pl = PL;
end
% End: Compute independent pdf  find

%%%%%%%%%%%%%%%%%%%%%%%%%%%%%%%%%%%%%%%%%%%%%%%%%%%%%%
% Check regularity and find wave characteristics and modify if  L3 ne 0
f.moments = [L0 L1 L2 L3 L4];
if isequal(L3,0)
    f.char = [L2/sqrt(L0*L4) sqrt(1-L2^2/(L0*L4))];
end
if modify && ~isequal(L3,0)
    fm = f;
    f.f(1,2:end) = fm.f(2,2:end);
    f.f(2:end,1) = fm.f(2:end,2);
    if nargout == 2
        fm = find;
        find.f(1,2:end) = fm.f(2,2:end);
        find.f(2:end,1) = fm.f(2:end,2);
    end
end
% End: Irregularity

return
\end{lstlisting}
}
\end{document}